\documentclass[a4paper,
               10pt,
               pdftex,
               normalheadings,
               headsepline,
               footsepline,
               headinclude,
               footinclude,
               DIV14,
               abstracton]
{scrartcl}

\usepackage[left=2.50cm, right=2.50cm, top=2.00cm, bottom=3.50cm]{geometry}

\usepackage
{
    graphicx,
    amssymb,
    amsmath,
    amsthm,
    xcolor,
    dsfont,
    algpseudocode,
    authblk,
}
\usepackage[ngerman, english]{babel}

\usepackage{enumerate}
\usepackage{algorithm}
\usepackage{tabularx}
\usepackage{subcaption} 
\usepackage{booktabs} 

\usepackage[colorlinks,
            pdffitwindow=false,
            plainpages=false,
            pdfpagelabels=true,
            pdfpagemode=UseOutlines,
            pdfpagelayout=SinglePage,
            bookmarks=false,
            colorlinks=true,
            hyperfootnotes=false,
            linkcolor=blue,
            citecolor=green!50!black]
{hyperref}

\usepackage{sectsty}
\sectionfont{\fontsize{10}{10}\selectfont}
\subsectionfont{\fontsize{10}{10}\selectfont}

\newcommand{\R}{\mathbb{R}}

\newcommand{\norm}[1]{\left\lVert#1\right\rVert}
\newcommand{\abs}[1]{\left\vert#1\right\vert}

\newcommand{\Conv}{\textnormal{{Conv}}}

\newcommand{\dt}{\,{d}t}
\newcommand{\dtau}{\,{d}\tau}
\newcommand{\Ev}{E_{\mathsf{V}}(V)}
\newcommand{\Emi}{E_{\mathsf{MI}}(\Delta t)}
\newcommand{\Emir}{E_{\mathsf{MI}^r}(\Delta t)}
\newcommand{\Emodel}{E_{\mathsf{model}}}
\newcommand{\Er}{E_{\mathsf{r}}(E)}
\newcommand{\Y}{Y}
\newcommand{\Fcal}{\mathcal{F}}
\newcommand{\Kcal}{\mathcal{K}}
\newcommand{\Lcal}{\mathcal{L}}

\newtheorem{theorem}{Theorem}[section]
\newtheorem{corollary}[theorem]{Corollary}
\newtheorem{lemma}[theorem]{Lemma}

\newtheorem{remark}[theorem]{Remark}
\newtheorem{example}[theorem]{Example}

\let\OLDthebibliography\thebibliography
\renewcommand\thebibliography[1]{
	\OLDthebibliography{#1}
	\setlength{\parskip}{0pt}
	\setlength{\itemsep}{0pt plus 0.3ex}
}

\usepackage{xcolor}
\newcommand{\rev}[1]{\textcolor{black}{#1}}
\newcommand{\revKaty}[1]{\textcolor{black}{#1}}

\begin{document}

\title{On the Universal Transformation of Data-Driven Models to Control Systems}
\author[1]{Sebastian Peitz}
\author[2]{Katharina Bieker}
																																			 
\affil[1]{\normalsize Department of Computer Science, Paderborn University, Germany}															 
\affil[2]{\normalsize Department of Mathematics, Paderborn University, Germany}

\date{}

\maketitle

\begin{abstract}
The advances in data science and machine learning have resulted in significant improvements regarding the modeling and simulation of nonlinear dynamical systems. It is nowadays possible to make accurate predictions of complex systems such as the weather, disease models or the stock market. Predictive methods are often advertised to be useful for control, but the specifics are frequently left unanswered due to the higher system complexity, the requirement of larger data sets and an increased modeling effort. In other words, surrogate modeling for autonomous systems is much easier than for control systems. In this paper we present the framework QuaSiModO (Quantization-Simulation-Modeling-Optimization) to transform arbitrary predictive models into control systems and thus render the tremendous advances in data-driven surrogate modeling accessible for control. Our main contribution is that we trade control efficiency by autonomizing the dynamics -- which yields mixed-integer control problems -- to gain access to arbitrary, ready-to-use autonomous surrogate modeling techniques. We then recover the complexity of the original problem by leveraging recent results from mixed-integer optimization. The advantages of QuaSiModO are a linear increase in data requirements with respect to the control dimension, performance guarantees that rely exclusively on the accuracy of the predictive model in use, and little prior knowledge requirements in control theory to solve complex control problems. 
\end{abstract}

Many challenges of our modern society could be addressed by significantly improving the control performance for highly complex systems in real-time, one example being the ever increasing demand for energy and the associated question of how it can be met in a cost-efficient and sustainable manner. The World Energy Outlook 2018 \cite{IEA18} suggests that -- assuming current policies remain unchanged and that the efficiency can be increased as expected -- the total demand in energy will increase by 25\% until 2040. Consequently, a further increase in efficiency of conventional processes for producing electrical energy (e.g., turbomachinery, combustion, wind energy and tidal energy) is required as well as the development of new concepts -- nuclear fusion being a prominent and very promising example \cite{PW05}. Furthermore, the efficiency of complex systems such as aircraft needs to be increased in order to limit the energy requirements.
All the above-mentioned applications are governed by high-dimensional nonlinear dynamics, which exhibit complex multi-scale phenomena and are thus extremely difficult to control. Similar challenges arise in other areas such as the health sector, where complex dynamics govern our breathing \cite{CKH02}, the flow of blood within arteries \cite{CRD+11} or of cerebrospinal fluid within the brain \cite{MTD+18}, and -- more recently -- the dynamics of pandemics such as COVID-19 \cite{GBB+20}.

The efficient prediction of complex systems such as the ones mentioned above is often hindered by the fact that the system dynamics are either very expensive to simulate or even unknown. Researchers have been investigating ways to accelerate the solution by using data for decades, the \emph{Proper Orthogonal Decomposition (POD)} being an early and very prominent example \cite{Sir87}. Therein, the original system is projected onto a linear subspace spanned by modes that have been computed from simulation data of the full system state.
More recently, the major advances in data science and machine learning have lead to a plethora of new possibilities to overcome the shortcomings of projection-based methods like POD, such as the required knowledge of the system dynamics and the quickly growing model dimension with increasing complexity. Important examples are different artificial neural network architectures such as \emph{Long Short-Term Memory (LSTM)} Networks \cite{HS97,VBW+18} or \emph{Reservoir Computers / Echo State Networks (ESN)} \cite{JH04,PHG+18}, regression-based frameworks for the identification of nonlinear dynamics \cite{BPK16,RBPK17}, or numerical approximations of the \emph{Koopman operator} \cite{RMB+09,Sch10,KNP+20}, which describes the linear dynamics of observable functions. These methods facilitate the efficient simulation and prediction of high-dimensional spatio-temporal dynamics using measurement data, without requiring prior system knowledge. One important application is, among many others, the prediction of rare events \cite{WVKS18}, for instance in nuclear fusion reactors \cite{KST19}.

The large success of data-driven modeling has also attracted the attention of the control community, where many approaches have been presented over the past decades. They can be categorized into the direct learning of feedback signals using, e.g., feed-forward neural networks \cite{SI18} or reinforcement learning \cite{DR11,LVV12,DFR15}, and into \emph{model predictive control (MPC)} methods \cite{GP17} using either intrusive approaches (i.e., the equations need to be known) such as POD \cite{KV99} or black-box methods, for instance specific neural network architectures \cite{BPB+20}. MPC allows for a particularly easy implementation of both control and state constraints.
However, a drawback is that the construction of surrogate models with inputs is often very tedious and in addition highly problem-specific. 
For instance, in the case of POD, multiple modeling steps and several simulations with different boundary conditions are required \cite{BCB05}. 
In pure black-box methods, a straight-forward approach to avoid this tedious effort is to transform control inputs via state augmentation \cite{PBK15}, by which the system is \emph{autonomized}, and then use ``off-the-shelf'' methods for autonomous systems. However, this results in significantly increased data requirements due to (a) the increased dimension of the augmented state space and (b) the fact that the dynamics of the control system is not necessarily restricted to a low-dimensional manifold any longer, which is a key enabler for efficient data-driven modeling \cite{BK19}.
Recently, an alternative was presented in \cite{PK19}, where the key observation was that instead of adapting the surrogate model according to the control problem requirements, it can be advantageous to modify the control problem and use a finite set of autonomous systems.
This results in a mixed-integer optimal control problem which is considerably harder to solve, and additional favorable properties of the surrogate model (such as the linearity of the Koopman operator \cite{POR20}) are required to facilitate an efficient solution.

\begin{figure}
	\centering
	\includegraphics[width=.7\textwidth]{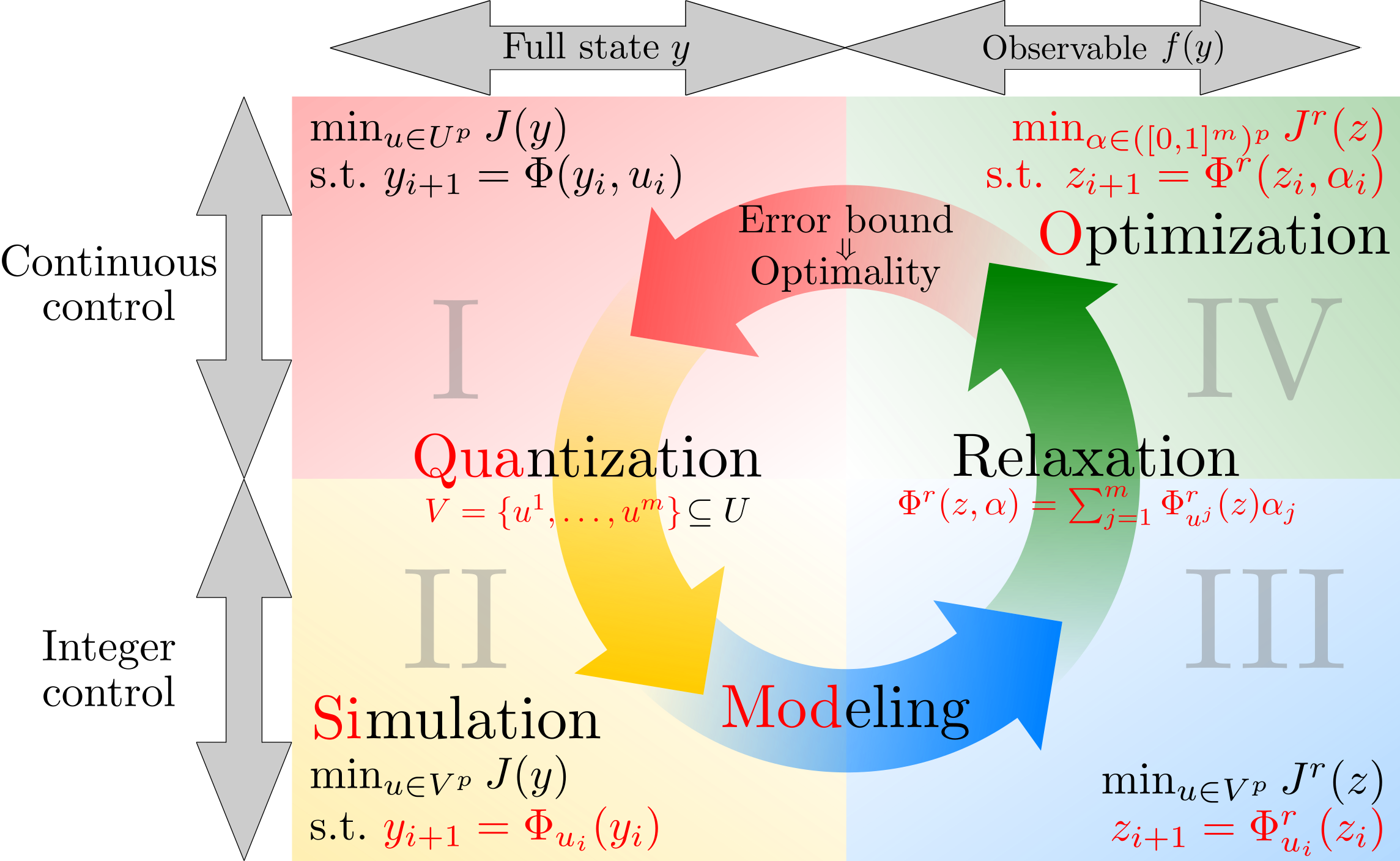}
	\caption{The QuaSiModO framework consisting of the four steps Quantization, Simulation, Modeling and Optimization.}
	\label{fig:QuaSiModO}
\end{figure}

The framework we present in this article -- \emph{QuaSiModO} -- makes use of the above-mentioned benefits of modifying the control problem. It consists of four main steps (cf.~also Figure \ref{fig:QuaSiModO}). In these, Problem \eqref{eq:OCP-I} (which is an optimal control problem with continuous inputs) is successively transformed into related control problems that -- as long as the predictive surrogate model is sufficiently accurate -- yield optimal trajectories $y^*$ that are close to one another. 
\begin{enumerate}[(i)]
    \item \underline{\textbf{Qua}ntization} of the the admissible control $U$ (for instance by replacing the interval $U=[u^{\mathsf{min}}, u^{\mathsf{max}}]$ by the bounds $V=\{u^{\mathsf{min}}, u^{\mathsf{max}}\}$);
    \item \underline{\textbf{Si}mulation} of the individual autonomous systems (e.g., $\Phi_{u^{\mathsf{min}}}(y) = \Phi(y,{u^{\mathsf{min}}})$ and $\Phi_{u^{\mathsf{max}}}(y) = \Phi(y,{u^{\mathsf{max}}})$); We may also use experimental data;
    \item \underline{\textbf{Mod}eling} of the individual systems -- using either the full state $y$ or some observable $z=f(y)$ -- via an arbitrary ``off-the-shelf'' surrogate modeling technique;
    \item \underline{\textbf{O}ptimization} using the set of autonomous surrogate models and relaxation techniques.
\end{enumerate}
This complex interplay between continuous and integer control modeling as well as between the full system state and observed quantities (e.g., measurements) allows us to utilize the best of both worlds:
\begin{itemize}
	\item integer controls for efficient data-driven modeling using arbitrary predictive models,
	\item continuous control inputs for real-time control, and
	\item existing error bounds for predictive models.
\end{itemize}

After introducing the respective steps in detail in the next section, we give a detailed derivation as well as numerical verification of the derived performance guarantees in Section \ref{sec:PerformanceGuarantees}. We then present control results on a large variety of complex systems, surrogate modeling techniques and observable types in Section \ref{sec:Results}.

\section{QuaSiModO steps}

The overarching goal of QuaSiModO is real-time control of complex systems using ``off-the shelf'' data-driven surrogate modeling techniques and their respective error bounds.
In the problems we consider, the state of such a complex dynamical system, denoted by $y$, is a function of time $t$ and -- in the case of partial differential equations (PDEs) -- of space $x$. For ordinary differential equations (ODEs) and PDEs, the dynamics is described by the right-hand side $g$, e.g., $\dot{y} = g(y(t), u(t))$ in the ODE case, with $u$ being the \emph{control input}. As most surrogate modeling techniques yield discrete-time systems, we directly introduce a time discretization using the time-T-map $\Phi$ with constant time step $\Delta t = t_{i+1} - t_i$, $i=0,1,\ldots$, i.e.,
\[
\Phi(y_i, u_i) = y_i + \int_{t_i}^{t_{i+1}} g(y(t), u_i) \dt = y_{i+1},
\]
where $u(t) = u_i\in U$ is constant over the interval $[t_i, t_{i+1})$, and $U$ is the admissible set for the control. 
A popular example for $U$ are box constraints: $U=[u^{\mathsf{min}}, u^{\mathsf{max}}]$. 
In the case of PDEs, we use a spatial discretization (such as finite elements), which then yields high-dimensional ODEs. Thus, we will only consider ODEs from now on.
Using the above considerations, the overall goal of MPC can be formalized in an optimal control problem of the following form:
\begin{equation} \tag{I} \label{eq:OCP-I}
\begin{aligned}
    \min_{u\in U^p} J(y) &= \min_{u\in U^p} \sum_{i=0}^{p-1} P(y_{i+1}) \\
    \mbox{s.t.} \quad y_{i+1} &= \Phi(y_i, u_i), \qquad i = 0,1,2,\ldots,
\end{aligned}
\end{equation}
where $J$ is the overall objective function over the time horizon $T = p\Delta t$, and $P$ is the objective function at a particular time instant, e.g., a tracking term $P(y_i) = \|y_i - y_i^{\mathsf{ref}}\|^2$. 
Note that no explicit dependency of $P$ on $u$ is assumed, as this will be favorable for the error analysis.
Furthermore, we assume that there exists a \emph{control-to-state} mapping that uniquely defines $y$ when choosing $u$, such that the state itself is not a control variable.

\subsection{Quantization:}
In the first step of QuaSiModO, we quantize the control set $U$ such that only a finite subset $V = \{u^1, \ldots, u^m \} \subseteq U$ is feasible. This allows us to replace the control system $\Phi(y,u)$ by a finite set of autonomous systems $\Phi_{u^j}(y)$, which yields the following mixed-integer problem:
\begin{equation} \tag{II} \label{eq:OCP-II}
\begin{aligned}
\min_{u\in V^p} J(y) &= \min_{u\in V^p} \sum_{i=0}^{p-1} P(y_{i+1}) \\
\mbox{s.t.} \quad y_{i+1} &= \Phi_{u_i}(y_i), \qquad i = 0,1,2,\ldots.
\end{aligned}
\end{equation}

\begin{remark}
The choice of the entries of $V$ is problem-specific and has an influence on the control performance and on the numerical effort. In addition, as we will see later, the reachable set corresponding to $U$ has to be contained in the convex hull of the reachable set corresponding to $V$ if we want to guarantee 
\revKaty{that the optimal value of $J$ in Problem~\eqref{eq:OCP-II} is arbitrarily close to the optimal value of Problem~\eqref{eq:OCP-I}}, cf.\ Lemma~\ref{lemma:existence_RCCP_solution}.
In the control-affine case, this is ensured by simply choosing the control bounds.
\end{remark}

\subsection{Simulation \& Modeling:}
Problem \eqref{eq:OCP-II} allows for the straightforward introduction of surrogate models, as we can simply replace the individual autonomous systems $\Phi_{u^j}(y)$ by respective reduced systems $\Phi^r_{u^j}(z)$ (or -- in continuous time -- replace $g$ by $g^r$) that predict the dynamics of a reduced quantity $z = f(y)$. 
Here, $f$ is called an \emph{observable} that takes arbitrary measurements from the full state. These range from full-state observation (i.e., $z=y$) over partial observations or point-wise measurements to arbitrary nonlinear functions of the state. In many applications -- in particular in experimental setups -- the choice of $f$ is determined by which measurements are feasible or accessible. For instance, in the case of fluid dynamics, it is highly unrealistic to measure the full state or even take measurements from the interior of the flow domain. 

The quantization also simplifies the process of simulation (or, more generally, data collection). Here, we can collect data from one long time series $Z = [z_1, \ldots, z_N]$ with random actuation $u_1,\ldots, u_N \in V$ (which we then have to split into data sets for the individual systems).

The fact that we need to construct surrogate models for a set of autonomous systems allows us to use arbitrary predictive models for $z$, thus giving us straightforward access to a highly active area of research with new techniques being proposed very frequently. These surrogate models then simply replace the full system in Problem~\eqref{eq:OCP-II}, yielding:
\begin{equation} \tag{III} \label{eq:OCP-III}
\begin{aligned}
    \min_{u\in V^p} J^r(z) &= \min_{u\in V^p} \sum_{i=0}^{p-1} P^r(z_{i+1}) \\
    \mbox{s.t.} \quad z_{i+1} &= \Phi_{u_i}^r(z_i), \qquad i = 0,1,2,\ldots.
\end{aligned}
\end{equation}
Note that the objective function has to be altered as well, as it is now a function of $z$ and not of the full state $y$. However, this is not a strong limitation, as objective function evaluations are usually based on observable quantities or the observables $f$ are chosen such that the objective can be calculated from the observed quantities.
This is particularly the case for real systems, where the measurements cannot be chosen arbitrarily, but need to be technically realized via sensors. Formally, we assume 
\begin{equation}\label{eq:EqualityObjectives}
P^r(f(y)) = P(y) \qquad \mbox{for all }y\in\Y.
\end{equation}

\subsubsection{Surrogate model error}
Regarding the \emph{quantification of the modeling error}, researchers have made significant advances over the past decades, starting with \emph{Proper Orthogonal Decomposition} \cite{Sir87}. 
With the increasing interest in data-driven methods, the effort invested in deriving bounds for entirely data-based models has significantly increased in recent years. For instance, convergence of Extended Dynamic Mode Decomposition towards the Koopman operator in the infinite data limit was shown in \cite{KM18b}, bounds for the prediction error with finite training data were presented in \cite{LT20,KNH20,NPP+21,ZZ21},
and probabilistic error bounds for the approximation of linear systems via finite data in \cite{SMT+18}.
Seeing the ever increasing effort, it is likely that more progress will be made in the near future. One of the goals of QuaSiModO is to exploit these advances in a systematic way for control, by deriving error bounds that exclusively rely on existing bounds for autonomous dynamics.
To this end, we assume in the following that an error bound for the deviation over a time step $\Delta t$ of the following form is known:
\begin{equation*}
\begin{aligned}            \norm{f(\Phi_{u^j}(y_0))-\Phi_{u^j}^r(\bar{z}_0)}_{\infty} \le E(\norm{f(y_0) - \bar{z}_0},\Delta t),
\end{aligned}
\end{equation*}
where $y_0$ denotes the initial value of the full state of the real system and $\bar{z}_0$ denotes the value the surrogate model is initialized with. In the first time step, $\bar{z}_0 = f(y_0)$ does hold, but this is usually not valid in subsequent steps due to the model error.
Thus, we obtain the time-dependent model error $\Emodel$:
\begin{equation}\label{eq:ErrorModel}
\begin{aligned}
\Emodel(t_0) &= \norm{f(y_0) - \bar{z}_0}, \\
\Emodel(t_i) &\le E(\Emodel(t_{i-1}),\Delta t).
\end{aligned}
\end{equation}

\subsection{Optimization:}
We now have three opportunities to calculate the solution of the surrogate-based Problem \eqref{eq:OCP-III}.
First, we can directly solve it (for instance, using Dynamic Programming) which -- due to the combinatorial nature -- is only viable for a small number of control parameters.
For larger problems, it is much more advisable to use a relaxation approach, which yields another continuous problem that can be solved efficiently using methods from nonlinear constrained optimization \cite{NW06}. 
The second and third approach thus rely on relaxation, and we can either directly apply the obtained solution to the real system or -- as this is only viable for control affine systems -- we can use the \emph{sum up rounding} (SUR) algorithm from \cite{SBD12}.  
This way, a control corresponding to one of the quantized inputs is applied to the real system. On the one hand, this limits the control freedom (which -- as we show in Section \ref{sec:PerformanceGuarantees} -- is acceptable for a wide range of problems). 
On the other hand, the real system is exclusively actuated by inputs that are already contained in $V$, which allows for a very easy implementation of online learning of the individual surrogate models. This can be extremely beneficial for the control performance, cf., e.g., \cite{BPB+20}.

As a first step towards the relaxation of \eqref{eq:OCP-III}, we introduce a new (binary) control variable $\omega \in (\{0,1\}^m)^p$ that yields a formulation equivalent to \eqref{eq:OCP-III}:
\begin{equation} \tag{III-$\omega$} \label{eq:OCP-IIIb}
\begin{aligned}
\min_{\omega \in (\{0,1\}^m)^p} &J^r(z) = \min_{\omega\in (\{0,1\}^m)^p} \sum_{i=0}^{p-1} P^{r}(z_{i+1}) \\
\mbox{s.t.} \qquad z_{i+1} &= \Phi^r(z_i, \omega_i) = \sum_{j=1}^m \omega_{i,j} \Phi^r_{u^j}(z_{i}),\\ \mbox{and} \quad\sum_{j=1}^m \omega_{i,j} &= 1, \qquad i = 0,1,2,\ldots.
\end{aligned}
\end{equation}
The last condition ensures that exactly one of the right-hand sides is applied in each time step, i.e., $\sum_{j=1}^m \omega_{i,j} u^j = u_i \in V$ for $i=0,\ldots, p-1$.

In order to obtain a continuous optimization problem which can be solved more efficiently than \eqref{eq:OCP-III} or \eqref{eq:OCP-IIIb}, where the control inputs are discrete, Problem \eqref{eq:OCP-IIIb} can be relaxed by replacing $\omega_i \in \{0,1\}$ with ${\alpha_i \in [0,1]}$.
\begin{equation} \tag{IV} \label{eq:OCP-IV}
\begin{aligned}
\min_{\alpha\in ([0,1]^m)^p} J^r(z) &= \min_{\alpha\in ([0,1]^m)^p} \sum_{i=0}^{p-1} P^{r}(z_{i+1}) \\
\mbox{s.t.} \quad z_{i+1} &= \Phi^r(z_i, \alpha_i) = \sum_{j=1}^m \alpha_{i,j} \Phi^r_{u^j}(z_i),\\  \mbox{and} \quad
&\sum_{j=1}^m \alpha_{i,j} = 1, \qquad i = 0,1,2,\ldots.
\end{aligned}
\end{equation}
The main differences to Problem \eqref{eq:OCP-I} are the additional condition on $\alpha$ and the fact that the dimension of the control variable is changed by a factor of $ (m-1) / \mathsf{dim}(U)$ (the last entry of $\alpha_i$ can always be computed from the condition $\sum_{j=1}^m \alpha_{i,j} = 1$). To solve \eqref{eq:OCP-IV}, efficient methods from nonlinear optimization such as gradient descent or BFGS methods \cite{NW06} can be used.

The final step of the optimization is now to construct the control $u^*$ from the optimal solution $\alpha^*$. 
To this end, we proceed according to the SUR procedure for mixed integer problems as proposed in \cite{SBD12}, 
where $\alpha^*$ is transformed back to $\omega^*$ while taking past rounding decisions into account:
\begin{subequations} \label{eq:SUR}
	\begin{align}
	\hat{\omega}_{i,j} &= \sum_{k=0}^{{i}} \alpha^*_{k,j} - \sum_{k=0}^{i-1} \omega^*_{k,j} \\
	\omega^*_{i,j} &= \begin{cases}1, & \hat{\omega}_{i,j} \ge \hat{\omega}_{i,l} ~\forall~ l \neq j \mbox{ and} \\ ~ &  j < l ~\forall~ l: \hat{\omega}_{i,j} = \hat{\omega}_{i,l}, \\
	0, & \mbox{else}.\end{cases} \\
	u_i^* &= \sum_{j=1}^m \omega^*_{i,j} u^j, \qquad i = 0,1,2,\ldots.
	\end{align}
\end{subequations}
For control-affine systems and convex sets $U$, an alternative, straightforward approach is to simply calculate the convex combination of the individual $u^j\in V$:
\[
u_i^* = \sum_{j=1}^m \alpha^*_{i,j} u^j, \qquad i = 0,1,2,\ldots.
\]
%
In the control-affine case, the two solutions are very close to one another, as experiments for the Lorenz system as well as the cylinder flow show, cf.\ Figure~\ref{fig:Results} in Section~\ref{sec:Results}.
\rev{Note that in terms of the objective function value, we can prove convergence for both the rounded and non-rounded approaches. However, an error bound on the control signal $u$ is not achievable in the rounded case.}

For the numerical solution of the resulting optimization problems, we simply use the standard optimizer from the \emph{SciPy} library, i.e., ``SLSQP'' \cite{Kra88} or ``trust-constr'' \cite{BHN99} with finite difference approximations of the derivatives. 
We note that significant speedups are very likely possible by using more efficient solvers tailored to the problems, as well as by providing explicit derivative information. 
However, the latter requires the derivative of the flow map $\Phi^r$ and thus, model-specific knowledge.
For neural networks, for instance, these gradients are often easily accessible via algorithmic differentiation.

\section{Performance guarantees}
\label{sec:PerformanceGuarantees}
In this section, we derive an error bound for the performance of the QuaSiModO approach in terms of the objective function value, which is composed of the quantization error, the model error and the error caused by the SUR procedure. Since the model error depends on the chosen model, we will assume that $\Emodel$ (cf.\ Eq.\ \eqref{eq:ErrorModel}) is known and concentrate on the quantization error which needs to be combined with the model error afterwards.
Note that even though we have introduced discrete-time systems, we will derive the bounds for continuous-time systems for ease of notation and in accordance with \cite{SBD12}. 
The discrete-time version follows as a special case.

\subsection{Similarity of trajectories}
We begin with the observation that the optimal solutions of \eqref{eq:OCP-I} and \eqref{eq:OCP-II} can get arbitrarily close. Obviously, the optimal value $J_{\eqref{eq:OCP-I}}^*$ of \eqref{eq:OCP-I} is always at least as small as the optimal value $J_{\eqref{eq:OCP-II}}^*$ of \eqref{eq:OCP-II}.
For the other direction, we will show that for each continuous control $u$, $u(t) \in U$, we can construct a sequence of discrete controls $v_n$, $v_n(t) \in V$, which leads to trajectories converging to the trajectory induced by $u$.
In \cite{W63} it was already stated that this holds if the set of points in the state space that can be reached with control inputs from $U$ is a subset of the convex hull of the points reachable with control inputs in $V$. We will prove the same, but in a constructive way in order to obtain concrete error bounds if this condition is not satisfied. Additionally -- and more importantly -- we are not in the limit case, i.e., a discretization error is introduced which depends on the switching time $\Delta t$ between the control inputs in \eqref{eq:OCP-II}. 

To estimate the distance between the trajectories given by different time-T-maps and different controls, the following result will be useful. A similar idea is used, for instance, in \cite{SBD12,Manns_PDE} to derive error bounds for the SUR algorithm.  
Unless otherwise specified, we denote by $\norm{\cdot}$ the maximum norm $\norm{\cdot}_{\infty}$ in the following, and $L^1((0,T),\Y)$ is the space of Lebesgue integrable functions mapping from the interval $(0,T)$ to the state space $Y$.

\begin{lemma}\label{lemma:distance_trajectories}
	Let $g,\bar{g}:\Y \times U \rightarrow \Y$ and $u,\bar{u}: [0,T] \rightarrow U$ be measurable functions with $\Y \subseteq \R^{n_y}$ and $U \subseteq \R^{n_u}$. 
	Furthermore, $g(y(\cdot),u(\cdot))$ and $\bar{g}(\bar{y}(\cdot), \bar{u}(\cdot)) \in L^1((0,T),\Y)$, where $y(\cdot)$ and $\bar{y}(\cdot)$ are given by
	\begin{equation*}
	\begin{aligned}
	y(t) &= y_0 + \int_0^t g(y(\tau),u(\tau)) \dtau\qquad\mbox{and}\\
	\bar{y}(t) &= \bar{y}_0 + \int_0^t \bar{g}(\bar{y}(\tau),\bar{u}(\tau)) \dtau,
	\end{aligned}
	\end{equation*}
	with $y_0, \bar{y}_0 \in \Y$. 
	Assume that $\bar{g}$ is Lipschitz continuous in the first argument with Lipschitz constant $L_{\bar{g}}$.
	If
	\begin{align*}
	\sup_{t \in [0,T]}\norm{\int_{0}^{t}{ g(y(\tau),u(\tau)) - \bar{g}(y(\tau),\bar{u}(\tau)) \dtau}} \le M
	\end{align*}
	then
	\begin{equation*}
	\norm{y(t)-\bar{y}(t)} \le (M +  \norm{y_0-\bar{y}_0}) e^{L_{\bar{g}} t} \quad \forall t \in [0,T].
	\end{equation*}
\end{lemma}
\begin{proof}
The claim mainly follows from Grönwall's lemma.\\
\end{proof}

In the MPC context, $y_0 = \bar{y}_0$ is the current state and we obtain a bound for the distance between two trajectories related to the two MPC problems \eqref{eq:OCP-I} and \eqref{eq:OCP-II}. Note that the given bound $M$ only has to hold for the optimal trajectory of the original MPC problem \eqref{eq:OCP-I}, not over the entire space. 

In order to prove an error bound between the optimal solutions of \eqref{eq:OCP-I} and \eqref{eq:OCP-II}, according to Lemma~\ref{lemma:distance_trajectories}, we need to construct a control $v$ with $v(t) \in V$ for the optimal solution $u^*$ of \eqref{eq:OCP-I} and derive the following bound:
\begin{align*}
\sup_{t \in [0,T]}\norm{\int_{0}^{t}{ g(y^*(\tau),u^*(\tau))- g(y^*(\tau),v(\tau))\dtau}} \le M.
\end{align*} 
To this end, we will use the idea of relaxation and SUR mentioned before, cf.\ Problems \eqref{eq:OCP-IIIb} and
\eqref{eq:OCP-IV}. 

\subsection{Bound without model error}
We begin without using a surrogate model, i.e., $z=f(y) = y$ and $\Phi^r = \Phi$, and
the error bound is derived in two steps. First, we prove that for every feasible solution of our original problem~\eqref{eq:OCP-I}, a feasible solution of the relaxed 
problem~\eqref{eq:OCP-IV} exists that leads to the same trajectory under some assumptions on the chosen subset~$V$. If these assumptions do not hold, we can give an error bound instead. Second, we use results from \cite{SBD12} to show that we obtain a trajectory of the 
binary control problem~\eqref{eq:OCP-IIIb} (which is equivalent to \eqref{eq:OCP-II} in this setting) that is arbitrarily close (for arbitrary small switching times) to the solution of~\eqref{eq:OCP-IV} by using SUR.

\begin{lemma}\label{lemma:existence_RCCP_solution}
	Let $U \subseteq \R^{n_u}$ be bounded and $V = \{u^1, \ldots, u^m \} \subseteq U$ a finite subset. Furthermore, let $g:\Y \times U \rightarrow \Y$ and $y: [0,T] \rightarrow \Y$ be  continuous and $u: [0,T] \rightarrow U$ be measurable.
	Then, there exists a measurable function $\alpha:[0,T]\rightarrow [0,1]^m$ such that $\sum_{j=1}^m \alpha_j(t) = 1~\forall t \in [0,T]$  and
	\begin{align*}
	\sup_{t\in [0,T]}&\norm{ \int_{0}^{t}g(y(\tau),u(\tau))-\sum_{j=1}^m g(y(\tau),u^j)\alpha_j(\tau)~\dtau} \\
	&\quad \le T \cdot  D =: M_1,
	\end{align*}
	where $D$ is the maximal distance between the reachable set corresponding to $U$ and the convex hull of the reachable set corresponding to $V$, more precisely, \[D = \sup_{t\in [0,T]} \sup_{\bar{y} \in g(y(t), U)} \operatorname{d}(\bar{y},\Conv(g(y(t),V)))\]
	with $\operatorname{d}(a,B) = \inf_{b \in B} \| a - b\|$ for $a \in \R^n$ and $B \subseteq \R^n$.
\end{lemma}

\begin{proof}
	For every $t\in [0,T]$, let  $y_c(t)$ be \rev{the} element in $\Conv(g(y(t),V))$ which is closest to $g(y(t),u(t))$, i.e., $\left\| g(y(t),u(t)) - y_c(t) \right\| = \operatorname{d}(g(y(t),u(t)),\Conv(g(y(t),V)))$.
	We can write $y_c(t)$ as a convex combination, i.e.,
	\begin{align*}
	y_c(t) =  \sum_{j=1}^m g(y(t),u^j) \alpha_j(t),
	\end{align*}
	with $\alpha(t) \in [0,1]^m$ and $\sum_{j=1}^m \alpha_j(t) = 1$. Note that, it is possible to choose $y_c(t)$ and $\alpha(t)$ such that $t \mapsto \alpha(t)$ is measurable. 
	Therefore, it holds for every $t \in [0,T]$
	\begin{align*}
	&\sup_{t \in [0,T]}\norm{\int_{0}^{t}g(y(\tau),u(\tau)) - \sum_{j=1}^m g(y(\tau),u^j)\alpha_j(\tau)~\dtau} \\
	\le & \enskip \sup_{t \in [0,T]}\int_{0}^{t}\norm{g(y(\tau),u(\tau))- y_c(\tau)} \dtau \\
	\le & \enskip T \cdot \sup_{t\in [0,T]} \norm{g(y(t),u(t))-y_c(t)}  \\
	\le & \enskip T \cdot \underbrace{\sup_{t\in [0,T]} \sup_{\bar{y} \in g(y(t), U)} \operatorname{d}(\bar{y},\Conv(g(y(t),V)))}_{ =~D~<~\infty, \text{ (continuity of $g$, $y$ and boundedness of } U\text{)}}= M_1.
	\end{align*}
\end{proof}

This means that if $V$ is chosen such that we can reach the extreme points of the reachable set corresponding to $U$ with the controls $u^j \in V$, we obtain $M_1 = 0$.
Next, we can use the idea of the SUR algorithm to estimate the error between the relaxed and the discrete control.
Therefore, we will use the following result.

\begin{theorem}\label{th:Sager}
	Let $g:\Y \times U \rightarrow \Y$, $y: [0,T] \rightarrow \Y$ and $\alpha:[0,T] \rightarrow [0,1]^m$ be measurable functions and assume that $\omega:[0,T] \rightarrow \{0,1\}^m$ is constructed from $\alpha$ via SUR (cf.\ Eq.\ \eqref{eq:SUR} or \cite{SBD12}) with time discretization $\Delta t$. Then,
	\begin{align*}
	\norm{\int_0^t \alpha(\tau)-\omega(\tau)} \le (m - 1) \Delta t\quad \forall t \in [0,T].
	\end{align*}
	Furthermore, assume that $g(y(\cdot),u^j)$ is differentiable for almost all $t \in [0, T]$ and that constants $C_1$ and $C_2 \in \R$ exist for all $u^j \in V$ such that for almost all $t \in [0, T]$:
	\begin{align*}
	\norm{\frac{d}{dt} g(y(t),u^j)} \le C_1\qquad\mbox{and}\qquad
	\norm{g(y(t),u^j)} \le C_2.
	\end{align*}
	Then,
	\begin{align*}
	\sup_{t \in [0,T]} &\norm{\int_{0}^{t}{\sum_{j=1}^mg(y(\tau),u^j)(\alpha_j(\tau)-\omega_j(\tau))}\dtau}\\
	& \le (C_2 +C_1T)(m - 1) \Delta t =:M_2(\Delta t).
	\end{align*}
\end{theorem}

\begin{proof}
	The proof can be found in \cite{SBD12}.
\end{proof}

\begin{remark}
	A similar result can be found in \cite{Manns_PDE} but with weaker assumptions on $g$. There, the authors prove that for $g(y(\cdot),v) \in L^1((0,T),\Y)$,
	\begin{align*}
	\lim_{\Delta t \rightarrow 0}&\sup_{t \in [0,T]}\norm{\int_{0}^{t}{\sum_{j=1}^m g(y(\tau),u^j)(\alpha_j(t)-\omega_j(\tau))}\dtau} =0\\ 
	&\mbox{if} \qquad
	\lim_{\Delta t \rightarrow 0} \norm{\int_0^t \alpha(\tau)-\omega(\tau)\dtau} = 0.
	\end{align*}
	Furthermore, they show that the statement holds not only for ODEs, but for semilinear PDEs as well.
	Here, we will use the result from \cite{SBD12}, 
	as we want to use the concrete error bound.
\end{remark}

Using the above results, we can now ensure that continuous and discrete inputs yield trajectories that are close.

\begin{theorem}\label{th:M_for_I_II}
	Let $U \subseteq \R^{n_u}$ be bounded and $V = \{u^1, \ldots, u^m \} \subseteq U$ be a finite subset.
	Assume $g:\Y \times U \rightarrow \Y$ is continuous, $u: [0,T] \rightarrow U$ is measurable and $y: [0,T] \rightarrow \Y$ is defined by
	\begin{align*}
	y(t) = y_0 + \int_0^t g(y(\tau),u(\tau)) \dtau, \quad y_0 \in \Y.
	\end{align*}
	Furthermore, let $g(y(t), u^j)$ be differentiable with respect to time for almost all $t \in [0, T]$ and assume that constants $C_1$ and $C_2 \in \R$ exist for all $u^j \in V$ such that for almost all $t \in [0, T]$:
	\begin{align*}
	\norm{\frac{d}{dt} g(y(t),u^j)} \le C_1\qquad\mbox{and}\qquad
	\norm{g(y(t),u^j)} \le C_2.
	\end{align*}
	In addition, assume that $g$ is Lipschitz continuous in the first argument with Lipschitz constant $L_g$.
	Then, for every $\varepsilon > 0$, there exists a discrete control function $\bar{u}: [0,T] \rightarrow V$, such that for  $\bar{y}$ given by
	\begin{align*}
	\bar{y}(t) = \bar{y}_0 + \int_0^t g(\bar{y}(\tau),\bar{u}(\tau)) \dtau, \quad \bar{y}_0 \in \Y, 
	\end{align*}
	it holds
	\begin{align*}
	\norm{y(t)-\bar{y}(t)} \le (M_1 + \varepsilon + \norm{y_0-\bar{y}_0}) \cdot e^{L_g t} \quad \forall t \in [0,T].
	\end{align*}
\end{theorem}

\begin{proof}
	Lemma~\ref{lemma:existence_RCCP_solution} ensures the existence of a measurable function $\alpha:[0,T] \rightarrow [0,1]^m$ with $\sum_{j=1}^m \alpha_j(t) = 1~\forall t \in [0,T]$ and $M_1 > 0$, such that
	\begin{align*}
	\sup_{t \in [0,T]}\norm{\int_{0}^{t}{g(y(\tau),u(\tau))-\sum_{j=1}^m g(y(\tau),u^j)\alpha_j(\tau)}\dtau} \le M_1.
	\end{align*}
	Using SUR, we can construct (according to Theorem~\ref{th:Sager}) a function $\omega:[0,T] \rightarrow \{0,1\}^m$ with\\ $\sum_{j=1}^m \omega_j(t) = 1~\forall t \in [0,T]$ from $\alpha$, such that
	\begin{align*}
	\sup_{t \in [0,T]} \norm{\int_{0}^{t}{\sum_{j=1}^m g(y(\tau),u^j)(\alpha_j(\tau)-\omega_j(\tau))}\dtau} \le M_2(\Delta t),
	\end{align*}
	with $M_2(\Delta t) = (C_2 +C_1T)(m - 1) \Delta t$.
	Therefore, we get 
	\begin{align*}
	\sup_{t \in [0,T]} &\norm{\int_{0}^{t}{ g(y(\tau),u(\tau))
			-\sum_{j=1}^m g(y(\tau),u^j) \omega_j(\tau)}\dtau}\\
		 &\le M_1 + M_2(\Delta t).
	\end{align*}
	Since $\bar{u}(t) := \sum_{j=1}^m \omega_j(\revKaty{t}) u^j \in V$ for all $t \in [0,T]$, choosing $\Delta t$ sufficiently small yields:
	\begin{align*}
	\sup_{t \in [0,T]} \norm{\int_{0}^{t}{ g(y(\tau),u(\tau))
			- g(y(\tau),\bar{u}(t))}\dtau} \le M_1 + \varepsilon,
	\end{align*} 
	and using Lemma~\ref{lemma:distance_trajectories}, we obtain the desired result.
\end{proof}

\begin{remark}
	The $m$ in $M_2(\Delta t)$ in Theorem~\ref{th:Sager} and \ref{th:M_for_I_II} can be reduced to the number of elements in $V$ that are actually required in the convex combinations to represent/approximate $u(t)$ over the prediction horizon $T$.
\end{remark}

Finally, with respect to the control problems \eqref{eq:OCP-I} and \eqref{eq:OCP-II}, we can derive a relation between the optimal values of Problems \eqref{eq:OCP-I} and \eqref{eq:OCP-II}.

\begin{corollary}
	Let $(y^*,u^*)$ be an optimal solution of \eqref{eq:OCP-I} where $U \subseteq \R^{n_u}$ is bounded and $P:\Y \rightarrow \R$ Lipschitz continuous with Lipschitz constant $L_P$ for a fixed initial value $y_0 \in \Y$.
	Assume $g:\Y \times U \rightarrow \Y$ and $V \subseteq U$ are of a form such that the requirements of Theorem~\ref{th:M_for_I_II} are satisfied. 
	Then, there exists a tuple $(\bar{y},\bar{u})$ with $\bar{u}(t) \in V$ which is feasible for \eqref{eq:OCP-II} with the same initial value $y_0$ such that
	\begin{align*}
	\abs{J(y^*)-J(\bar{y})} \le L_P (M_1 +  M_2(\Delta t)) \frac{e^{L_g \Delta t} \left(e^{p L_g \Delta t}-1\right)}{e^{L_g \Delta t}-1}.
	\end{align*}
\end{corollary}

\begin{proof}
The claim follows from Theorem~\ref{th:M_for_I_II} and the Lipschitz continuity of $P$.
\end{proof}

\begin{remark}
	The errors 
	\begin{equation}\label{eq:ErrorsEV_Emi}
	\begin{aligned}
	\Ev&=L_P M_1 \frac{e^{L_g \Delta t} \left(e^{p L_g \Delta t}-1\right)}{e^{L_g \Delta t}-1}\\ 
	\mbox{and}\quad \Emi&=L_P M_2(\Delta t) \frac{e^{L_g \Delta t} \left(e^{p L_g \Delta t}-1\right)}{e^{L_g \Delta t}-1}
	\end{aligned}
	\end{equation}
	account for the distance between the reachable sets of $V$ and $U$ and the transformation into a mixed integer problem and the corresponding relaxation technique, respectively. For an appropriate choice of $V$, we have $M_1=0$, cf.\ Lemma~\ref{lemma:existence_RCCP_solution}. Moreover, we have
	\[
	\lim_{\Delta t \rightarrow 0}\frac{e^{L_g \Delta t} \left(e^{p L_g \Delta t}-1\right)}{e^{L_g \Delta t}-1} = p \quad \mbox{and} \quad \lim_{\Delta t \rightarrow 0} M_2(\Delta t)=0,
	\]
	such that both errors can become arbitrarily small using the appropriate numerical setup.
\end{remark}
\begin{remark}
	An error bound can also be obtained if $P$ (and $J$) explicitly depend on the control $u$. In this case, however, an additional term needs to be added that pessimistically bounds the distance between the optimal controls for the full and surrogate-based problems. This bound also depends on the size of the control set $U$.
\end{remark}  

\subsection{Combination with model error}
We now consider the additional errors resulting from surrogate models (Eq.\ \eqref{eq:ErrorModel}), i.e., we solve Problem \eqref{eq:OCP-III} with the assumption that $J$ and $J^r$ are equivalent, cf.\ Eq.\ \eqref{eq:EqualityObjectives}. To do so, there are several ways in practice, and depending on the solution strategy, we obtain different error bounds. One can immediately see (via the triangle inequality) that the different error sources are additive in all cases.
For ease of notation, we introduce the \emph{control-to-state operators} $S: U^p \rightarrow \Y^p$ and $S^r: U^p \rightarrow f(\Y)^p$ which -- for a fixed $y_0 \in \Y$ -- map the control inputs $(u_0,\dots,u_{p-1})$ to the corresponding states $(y_1,\dots,y_p)$ and $(z_1,\dots,z_p)$, respectively. 

\subsubsection{Nonlinear observable functions}
The obvious approach is to solve problem \eqref{eq:OCP-III} directly, for instance using Dynamic Programming or a total evaluation of all possible combinations of control inputs. If we have an appropriate (i.e., sufficiently small) set $V$, this can be very efficient. 
In this case, we need to bound the difference  $\abs{J(S(u^*_{\eqref{eq:OCP-I}}))- J(S(u^*_{\eqref{eq:OCP-III}}))}$, where $u^*_{\eqref{eq:OCP-I}}$ and $u^*_{\eqref{eq:OCP-III}}$ are the optimal solutions of \eqref{eq:OCP-I} and \eqref{eq:OCP-III}, respectively. 
Since we already have a bound for $\abs{J(S(u^*_{\eqref{eq:OCP-I}}))- J(S(u^*_{\eqref{eq:OCP-II}}))}$, it is sufficient to determine the error $\abs{J(S(u^*_{\eqref{eq:OCP-II}}))- J(S(u^*_{\eqref{eq:OCP-III}}))}$. Therefore, we first consider the error between $J(S(u))$ and $J^r(S^r(u))$ for arbitrary $u \in V^p$:
\begin{align}
|J(S(u)) &- J^r(S^r(u))| = \abs{J^r(f(S(u))) - J^r(S^r(u))} \notag \\
&= \abs{\sum_{i=0}^{p-1} P^r(f((S(u))_i)) - P^r((S^r(u))_i)} \notag \\ 
&\le  \sum_{i=0}^{p-1} L_{P} \norm{f((S(u))_i) - (S^r(u))_i} \label{eq:model_error_I_III} \\
&= \sum_{i=0}^{p-1} L_{P} \norm{f\left( \Phi_{u_i}(y_{i})\right)- \Phi^r_{u_i}(z_{i})} \notag \\
&\le {L_P} \sum_{i=1}^{p}  \Emodel(t_i), \notag 
\end{align}
where the first equality is due to \eqref{eq:EqualityObjectives}.
For ease of notation, we do not distinguish between the Lipschitz constants of $P$ and $P^r$ and just use maximum of the two as $L_P$. 
Now, it holds
\newcommand{\aStar}{u^*_{\eqref{eq:OCP-II}}}
\newcommand{\aBar}{u^*_{\eqref{eq:OCP-III}}}
\begin{align*}
&|J(S(u^*_{\eqref{eq:OCP-III}})) - J(S(u^*_{\eqref{eq:OCP-II}}))| = J(S(u^*_{\eqref{eq:OCP-III}}))- J(S(u^*_{\eqref{eq:OCP-II}}))\\
&= \underbrace{J(S(\aBar)) - J^r(S^r(\aBar))}_{\le \abs{J(S(\aBar))) - J^r(S^r(\aBar))}} \\&+ \underbrace{J^r(S^r(\aStar)) - J(S(\aStar))}_{\le \abs{J^r(S^r(\aStar)) - J(S(\aStar)))}}\\
&+ \underbrace{(J^r(S^r(\aBar)) - J^r(S^r(\aStar))}_{\le 0\text{, since $\aBar$ is a global Min. of $J^r(S^r(\cdot))$}}\\
&\le ~2 {L_P} \sum_{i=0}^{p}  \Emodel(t_i),
\end{align*}
where we can neglect the absolute value because $u^*_{\eqref{eq:OCP-II}}$ is optimal with respect to $J$.
In summary, the error bound is thus given by 
\begin{equation}\label{eq:ErrorIII}\tag{E1}
\begin{aligned}
|&J(S(u^*_{\eqref{eq:OCP-I}}))- J(S(u^*_{\eqref{eq:OCP-III}}))| \\
&\le\Ev + \Emi + \underbrace{2 {L_P} \sum_{i=0}^{p} \Emodel(t_i)}_{\Er},
\end{aligned}
\end{equation}
see Eq.\ \eqref{eq:ErrorModel} for $\Emodel$ and Eq.\ \eqref{eq:ErrorsEV_Emi} for $\Ev$ and $\Emi$.

In many practical applications, we will have a finite set $V = \{u^1,\dots,u^m\} \subseteq U$ which is too large to solve the combinatorial problem \eqref{eq:OCP-III} directly. 
In this case, we can solve the relaxation \eqref{eq:OCP-IV} in combination with SUR instead, where the resulting control is denoted by $u^*_{\eqref{eq:OCP-IV}-\mathsf{SUR}}$. 
The error is composed of the error \eqref{eq:ErrorIII} and $\abs{J(S(u^*_{\eqref{eq:OCP-III}}))- J(S(u^*_{\eqref{eq:OCP-IV}-\mathsf{SUR}}))}$. Analog to the derivation before, we obtain an error $\Emir$ for the surrogate-based mixed-integer transformation (with variables $M_2^r(\Delta t)$ and $L_{g^r}$ if the surrogate model is given by $\Phi^r(y_i, u_i) = y_i + \int_{t_i}^{t_{i+1}} g^r(y(t), u_i) \dt$ and $g^r$ is Lipschitz continuous with Lipschitz constant $L_{g^r}$). Therefore, we get
\begin{equation}\label{eq:ErrorIIIbSUR}\tag{E2.a}
\begin{aligned}
|&J(S(u^*_{\eqref{eq:OCP-I}}))- J(S(u^*_{\eqref{eq:OCP-IV}-\mathsf{SUR}}))|\\
 &\le \Ev + \Emi + 2\Er + \Emir,
\end{aligned}
\end{equation}
where the additional second model error $\Er$ is due to the fact that we need to estimate the difference  in terms of the true objective.

\begin{table*}[h]
	\centering
	\caption{Different transformation procedures and corresponding error bounds.}
	\begin{tabular}{l|c|c|c|l}
		Transformation approach & Error bound & control & $f$ linear &  Type of \\
		& & affine & & optimization \\
		\hline
		\eqref{eq:OCP-I} $\rightarrow$ \eqref{eq:OCP-II} $\rightarrow$ \eqref{eq:OCP-III} & Eq.~\eqref{eq:ErrorIII} & --- & --- & Combinatorial \\
		\eqref{eq:OCP-I} $\rightarrow$ \eqref{eq:OCP-II} $\rightarrow$ \eqref{eq:OCP-III} $\rightarrow$ \eqref{eq:OCP-IV} $\rightarrow$ SUR & Eq.~\eqref{eq:ErrorIIIbSUR} & --- & --- & Continuous \\
		& Eq.~\eqref{eq:ErrorIIIbSUR_linOb} & --- & $\checkmark$ & Continuous \\
		\eqref{eq:OCP-I} $\rightarrow$ \eqref{eq:OCP-II} $\rightarrow$ \eqref{eq:OCP-III} $\rightarrow$ \eqref{eq:OCP-IV} & Eq.~\eqref{eq:ErrorIV} & $\checkmark$ & $\checkmark$ & Continuous
	\end{tabular}
	\label{tab:Transformations}
\end{table*}

\subsubsection{Linear observable functions}
If $f$ is linear (which is in particular the case when considering the full-state observable), an error bound can be obtained in a similar way to Eq.~\eqref{eq:ErrorIII}. Therefore, we do a similar computation as in \eqref{eq:model_error_I_III} for the relaxed systems.
In this case, we obtain the same bound as for \eqref{eq:ErrorIII}:
\begin{equation}\label{eq:ErrorIIIbSUR_linOb}\tag{E2.b}
\begin{aligned}
|&J(S(u^*_{\eqref{eq:OCP-I}}))- J(S(u^*_{\eqref{eq:OCP-IV}-\mathsf{SUR}}))| \\
&\le \Ev + \Emi + \Er.
\end{aligned}
\end{equation}

Finally, the third option is to solve \eqref{eq:OCP-IV} and directly apply the relaxed solution to the original system. Obviously, this is only feasible if $\Conv(V) \subseteq U$. This way, we introduce an additional error caused by the linear interpolation if the system is not control affine. Nevertheless, many systems are control affine and in this case, we have
\begin{align}\label{eq:ErrorIV}\tag{E3}
\abs{J(S(u^*_{\eqref{eq:OCP-I}})) - J\left(S\left(\sum_{j=1}^m \alpha^*_{\eqref{eq:OCP-IV},j} u^j\right)\right)} \le \Ev + \Er.
\end{align}

The bounds are summarized in Table \ref{tab:Transformations}, together with the additional requirements for the control problem. It should be noted that all errors besides the one for the surrogate model can be made arbitrarily small. We have $\Ev = 0$ if the convex hull of the reachable set corresponding to $V$ is a subset of the one corresponding to $U$, and both $\Emi$ and $\Emir$ vanish as the switching time $\Delta t$ tends to zero.

\subsection{Example}\label{subsec:ErrorboundExample}

To study the error bounds numerically, we consider the well-known Duffing oscillator:
\begin{align*}	  
\dot{y} = g^r(y,u)=\begin{pmatrix} y_2 \\ -\delta y_2 - \alpha y_1 - \beta y_1^3 + \varepsilon \end{pmatrix} + \begin{pmatrix} 0 \\ u \end{pmatrix},
\end{align*}
with constants $\alpha=-1$, $\beta=1$, $\delta=0$, and $u(t)\in U=[-4,4]$. To introduce a model error, we add a constant perturbation in the second equation, i.e., 
we set $\varepsilon = 0$ and $\varepsilon = 10^{-1}$ for the true and the surrogate model, respectively.
As the finite set of controls we choose $V=\{-4,4\}$. 
To determine the model error $\Emodel(t_i)$ in \eqref{eq:ErrorModel} we can use Lemma~\ref{lemma:distance_trajectories} since we use the full-state observable. For an arbitrary control $u(t)$ and a trajectory $y(t)$, it holds
\begin{align*}
\sup_{t \in [0,\Delta t]} \norm{\int_0^t g(y(\tau),u(\tau)) -  g^r(y(\tau),u(\tau))} \le \varepsilon \Delta t
\end{align*}
and according to Lemma~\ref{lemma:distance_trajectories}, we can estimate the model error via
\begin{align*}
\Emodel(t_{i+1}) &= \norm{y_{i+1} - y^r_{i+1}} \le (\varepsilon \Delta t +  \norm{y_{i} - y^r_{i}}) e^{L_g \Delta t}\\
&= ( \varepsilon \Delta t + \Emodel(t_{i}))e^{L_{g} \Delta t},
\end{align*}
where $y_i$ and $y^r_i$ are the discrete trajectories defined by $g$ and $g^r$, respectively, for a given starting point $y_0 \in \Y$ and a control sequence $u_i \in V$.  
The goal is to stabilize the system at $y^{\mathsf{ref}} = (0,0)^\top$, i.e., $J(y) = \sum_{i=0}^{p-1} \norm{y_{i+1}}$, and we solve Problems \eqref{eq:OCP-I}, \eqref{eq:OCP-IV} and \eqref{eq:OCP-III} (the latter using SUR) to investigate the error bounds \eqref{eq:ErrorIV} and \eqref{eq:ErrorIIIbSUR_linOb}. The control horizon is $[0,1]$, and we use $\Delta t = 2\cdot 10^{-3}$ for the time-T-maps $\Phi$ and $\Phi^r$, 
as well as for the SUR.
%

The constants $C_1$ and $C_2$ that enter the calculation of $M_2$ are approximated from data using several simulations with random initial conditions. We do the same for the Lipschitz constants, which we estimate via the derivative. As the system is control affine and $V$ consists of $u_{\mathsf{min}}$ and $u_{\mathsf{max}}$, we have $M_1=0$. 

\begin{figure}[t!]
	\centering
	 \includegraphics[width=0.5\columnwidth]{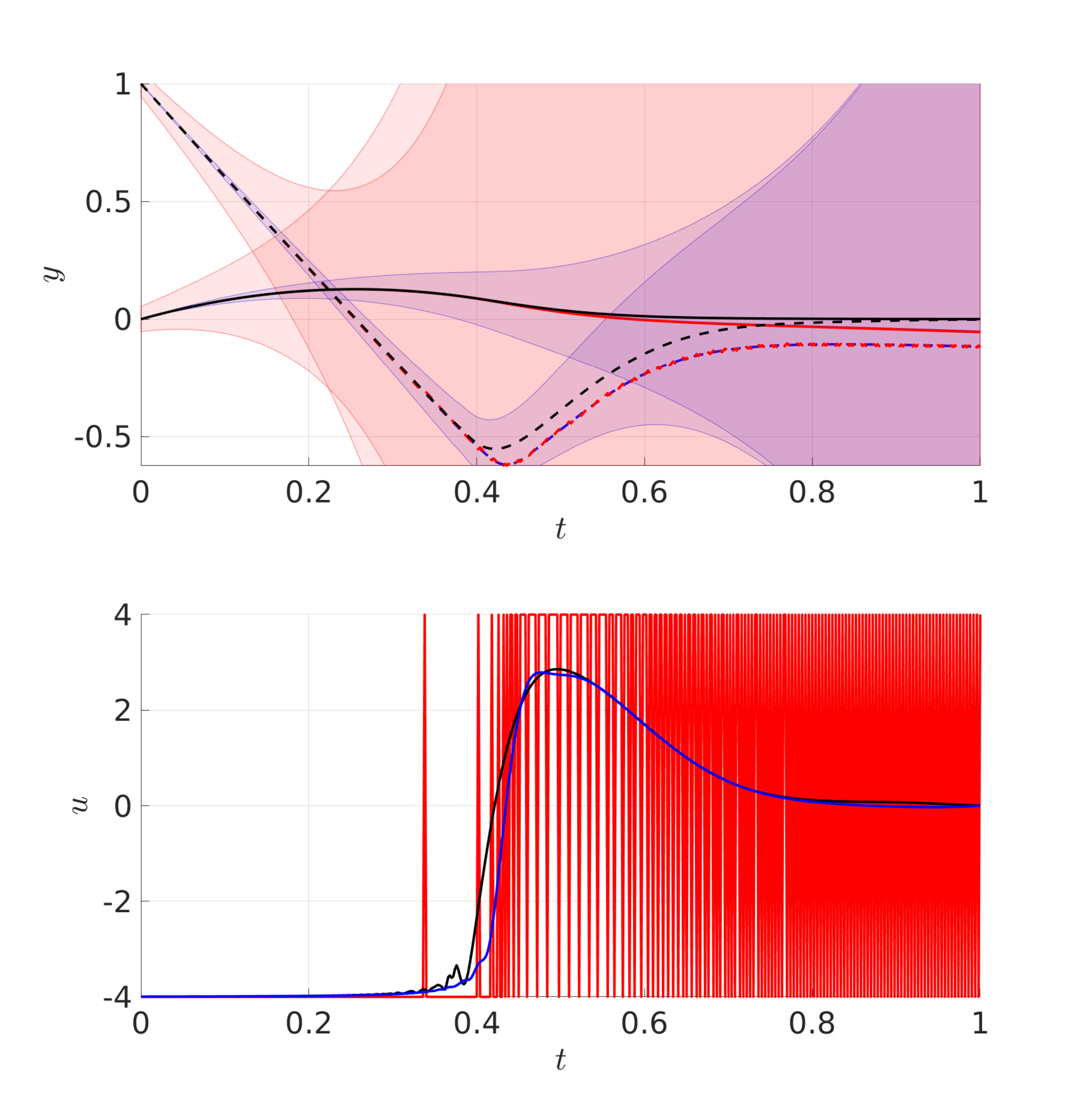} 
	 \caption{Top: Trajectory of the interpolated control (blue) and the SUR approach (red) with corresponding error bounds \eqref{eq:ErrorIV} and \revKaty{\eqref{eq:ErrorIIIbSUR_linOb}}. The solid lines denote the component $y_1$ and the dashed lines $y_2$. Bottom: Corresponding control inputs. The optimal state $y^*$ and control $u^*$ of Problem \eqref{eq:OCP-I} are shown in black.}
	\label{fig:Duffing_yu}
\end{figure}

Figure~\ref{fig:Duffing_yu} shows the results for the two approaches. We see that both achieve the control task relatively well with a small error due to the constant offset $\varepsilon$ in the second component. Moreover, we see that the error bound \eqref{eq:ErrorIV} is very well suited for the MPC context (where short prediction horizons are very common), and it is much tighter than the SUR approach, even though
the two errors come closer with decreasing $\Delta t$. Nevertheless, it can be concluded that solving Problem \eqref{eq:OCP-IV} without rounding is clearly advantageous for control-affine systems.

\begin{figure*}[h!]
	\centering
	\includegraphics[width=\textwidth]{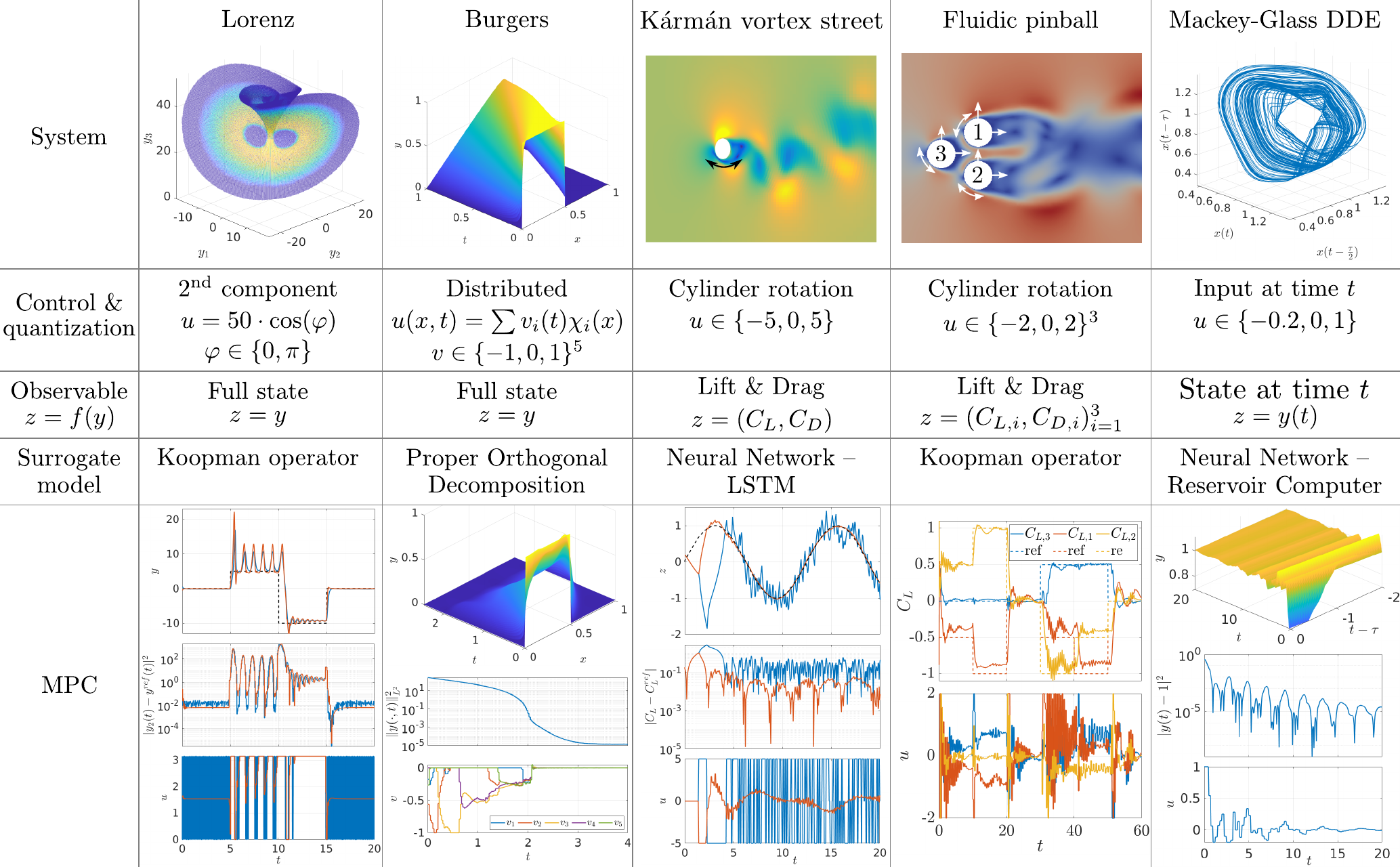}
	\caption{MPC using QuaSiModO applied to various combinations of systems and surrogate models. In the last row, the results of the MPC are shown for the settings explained above. For each system, the first row shows the state space, the second one the tracking error and the last one the obtained control over time. For the Burgers equation and the Mackey-Glass system, the control task was to stabilize the system at $0$ and $1$, respectively. For the other systems, the reference trajectory is shown by the dashed lines in the $y$ plot. In all examples, problem \eqref{eq:OCP-IV} was solved. For the Lorenz and Navier--Stokes example, both the continuous (orange) and the control computed via SUR (blue) are shown. In the remaining examples, we have applied the continuous solution of Problem \eqref{eq:OCP-IV} to the real system. Higher-dimensional inputs are considered for the Burgers equation ($n_u = 5$) and the fluidic pinball ($n_u=3$).}
	\label{fig:Results}
\end{figure*}

\section{Results}
\label{sec:Results}
We have tested the QuaSiModO framework on a variety of dynamical systems, observable functions and surrogate modeling techniques, cf.~Figure~\ref{fig:Results}, a detailed description of the numerical setup is given in Appendix \ref{app:modeling} for the surrogate modeling and in Appendix \ref{app:numerics} for the test problems. The systems range from the chaotic Lorenz system over the dynamics of blood cells (modeled by the Mackey-Glass delay differential equation) to the Navier--Stokes equations for fluid flows.
The observed quantities $z$ range from the full state over partial observations and delay coordinates to point-wise measurements, and we use surrogate models based on POD, the Koopman operator and different neural network architectures. 
For instance, we can control the lift force acting on a cylinder (determined by the velocity and pressure fields governed by the 2D Navier--Stokes equations) without any knowledge of the flow field using the standard LSTM framework included in \emph{TensorFlow},
and stabilize the Mackey-Glass equation using a standard ESN.
This highlights the flexibility and broad applicability of the method and the success of the technique in constructing data-driven feedback controllers. 

In the following, we describe \rev{three} experiments in more detail. First, the Lorenz system is analyzed to present the difference between affine linear and nonlinear control systems. \rev{We then briefly comment on the fluidic pinball, 
\revKaty{which has a higher-dimensional input and possesses the most complex dynamics from all our experiments.}
Finally,} the Mackey Glass equation is described to illustrate the potential to exploit certain surrogate model structures. For reconstruction of the results, we refer to the Code and its documentation which is available online, as well as the more detailed description in the appendices.

\subsection{Lorenz system \& Koopman operator:}\label{subsec:Lorenz}
\begin{figure}[h!]
\centering
\parbox[b]{0.45\textwidth}{
    \centering
	\includegraphics[width=0.45\textwidth]{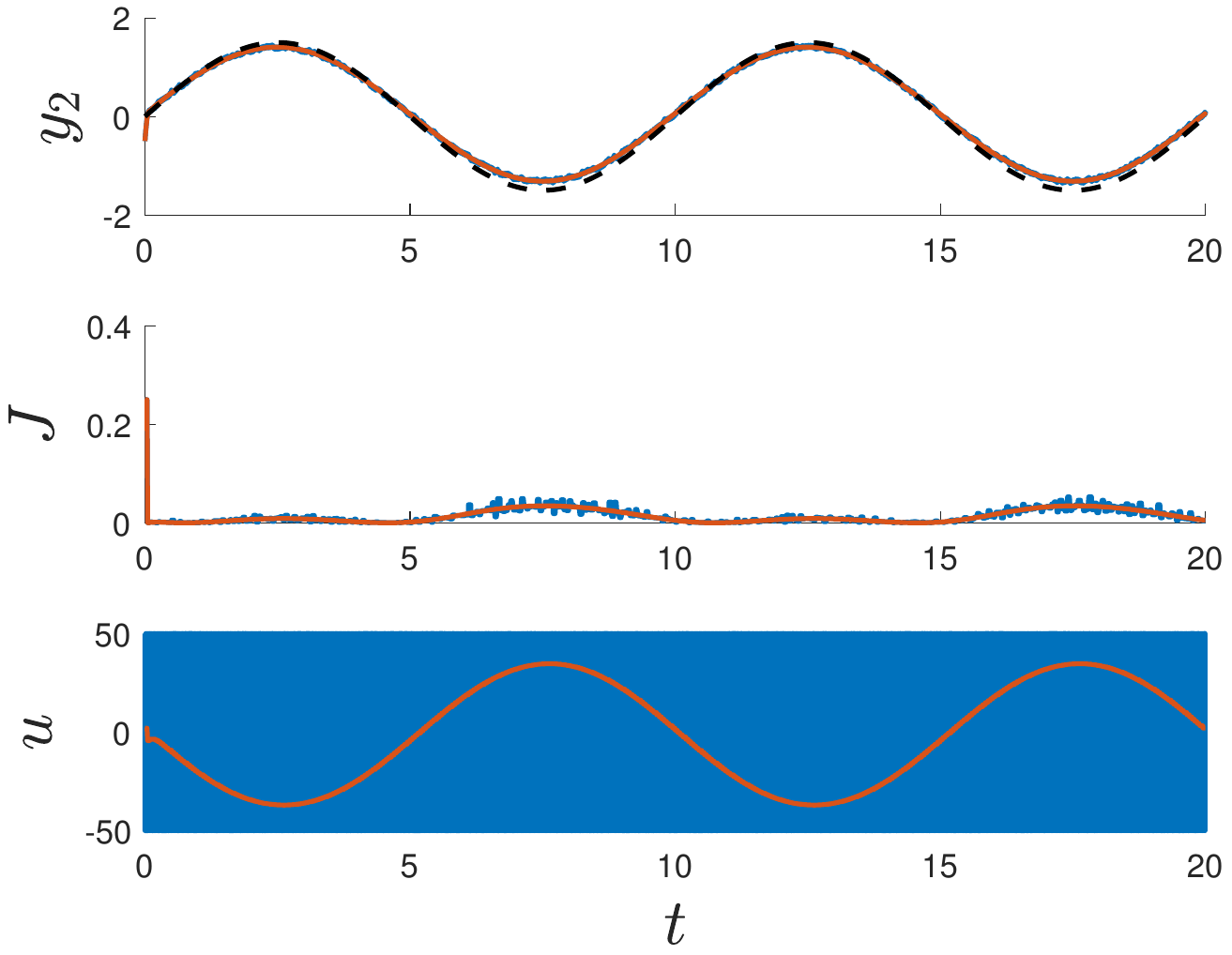}\\
	(a)
}
\hfill
\parbox[b]{0.45\textwidth}{
    \centering
	\includegraphics[width=0.45\textwidth]{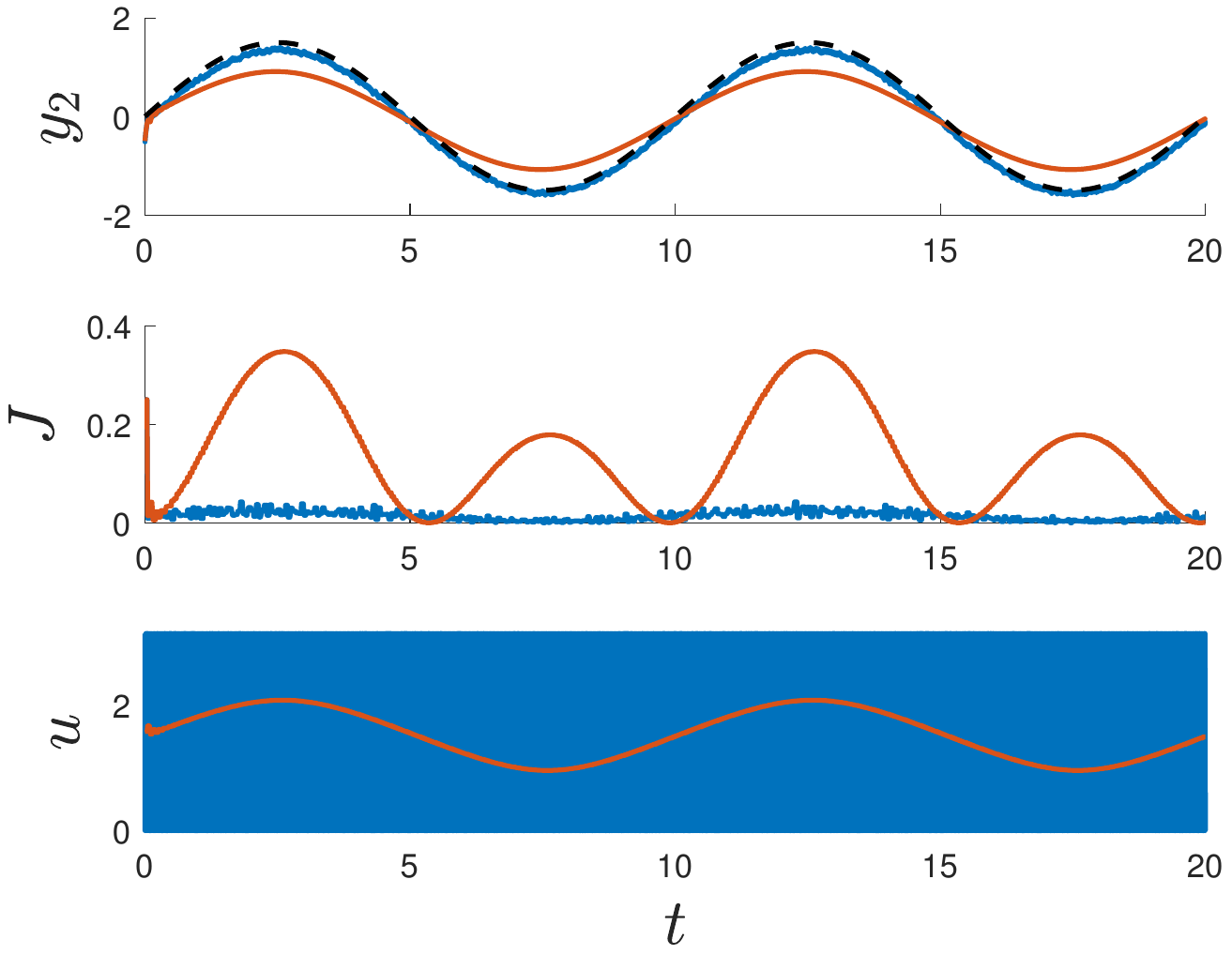}\\
	(b)
}
\caption{Solutions for the control affine (a) and the nonlinear Lorenz system (b) with an EDMD surrogate model. The solution of the linearly interpolated problem is shown in orange and the result generated via SUR after each optimization step in blue. The reference trajectory ($y_{2,i}^{\textnormal{ref}} = 1.5 \cdot \sin(4\pi \cdot t_i / T_{\mathsf{MPC}})$) is given by the dashed black line.} \label{fig:Lorenz_EDMD}
\end{figure}

As a first example, we consider the Lorenz system.
Despite the chaotic dynamics, there are many studies showing that predictions over time horizons of moderate length are possible using different methods such as neural networks \cite{CHS19} or sparse regression \cite{BPK16,KKB18}.
Here, the control task is to track a reference trajectory for the second variable, i.e., $P(y_i) = \norm{y_{2,i}^{\textnormal{ref}}-y_{2,i}}_2^2$, using an additive control input:
\begin{equation}\label{eq:Lorenz_linear}
\begin{aligned}
    \frac{d}{dt}\begin{pmatrix} y_1 \\ y_2 \\ y_3 \end{pmatrix} = \begin{pmatrix} \sigma (y_2 - y_1) \\ y_1 (\rho- y_3) - y_2 \\ y_1 y_2 - \beta y_3 \end{pmatrix} + \begin{pmatrix} 0 \\ u \\ 0 \end{pmatrix},
\end{aligned}
\end{equation}
where $(\sigma, \rho, \beta) = (10, 28, \frac{8}{3})$.
As the surrogate model, we use the Koopman operator approximated via EDMD \cite{WKR15,KKS16}, and the full state observable $z=f(y) = y$. 
The control sets are $U=[-50,50]$ and $V=\{-50,50\}$, respectively, and we collect 2000 data points over 100 seconds ($\Delta t=0.05$) using piecewise constant inputs from $V$.
As the control input enters linearly, the solution of \eqref{eq:OCP-IV} can be directly applied to the real system without an additional error. Therefore, it is not surprising that we achieve a very good control performance using both the interpolated and rounded solutions, cf. Figure~\ref{fig:Lorenz_EDMD} (a).

In order to emphasize the importance of the effect of nonlinear control inputs, we study a variation of the problem, i.e.,
\begin{equation}\label{eq:Lorenz_nonlinear}
\dot{y}_2 = y_1 (\rho - y_3) - y_2 + 50 \cdot \cos(u),
\end{equation}
with $U = [0, \pi]$, and chose $V = \{0, \pi\}$ to ensure $M_1 = 0$, cf.\ Lemma~\ref{lemma:existence_RCCP_solution}. Here, an interpolation error is introduced by applying the linearly interpolated control to the system and therefore, the SUR algorithm provides a substantially better solution, cf.~Figure~\ref{fig:Lorenz_EDMD} (b). 

\rev{
\subsection{Fluidic pinball \& Koopman operator:}\label{subsec:FP}
The \emph{fluidic pinball} is an unsteady flow around three cylinders arranged on an equidistant triangle, governed by the two-dimensional incompressible Navier--Stokes equations. 
In \cite{DNMP19}, it was conjectured that the system exhibits chaotic behavior above the critical Reynolds number $Re\approx 115$. In our control setting, we use $Re= 200$ as in \cite{BPB+20,POR20}, which likely results in the most complex problem from a dynamics perspective. The goal is to perform tracking of the lift coefficients of the three cylinders, and the control input is the rotation of all cylinders, meaning that we have a three-dimensional input. 
For the surrogate modeling, we proceed in the same way as for the Lorenz systen, and approximate the Koopman operator using EDMD. We do not observe the full state, but only the lift and drag of the three cylinders, which allows for a much more efficient (and more realistic, from an engineering perspective) modeling.
As the fluid solver, we use \emph{OpenFOAM}, and the setup is identical to the references above. Moreover, the example can be found in the accompanying git repository {https://github.com/SebastianPeitz/QuaSiModO}.}

\subsection{Mackey-Glass equation \& \revKaty{Echo State Network}}
\label{subsec:MG}
The second example is the control of the Mackey-Glass equation which is a delay differential equation modeling blood cell reproduction: 
\begin{equation}
    \begin{aligned}
    \dot{y}(t) = \beta \frac{y(t -\tau)}{1 + y(t -\tau)^\eta} - \gamma y(t) + u(t).
    \end{aligned}
\end{equation}
The additive term to control the system represents -- in the context of blood reproduction -- an increase in the number of blood cells caused by, e.g., a transfusion.
The uncontrolled system (i.e., $u(t) = 0$) was studied for different parameters and chaotic behavior was proven for certain parameter values, for instance $\beta = 2$, $\gamma = 1$ and $\eta = 9.65$ \cite{GM79}.
\begin{figure}
	\centering
	\includegraphics[width=0.6\textwidth]{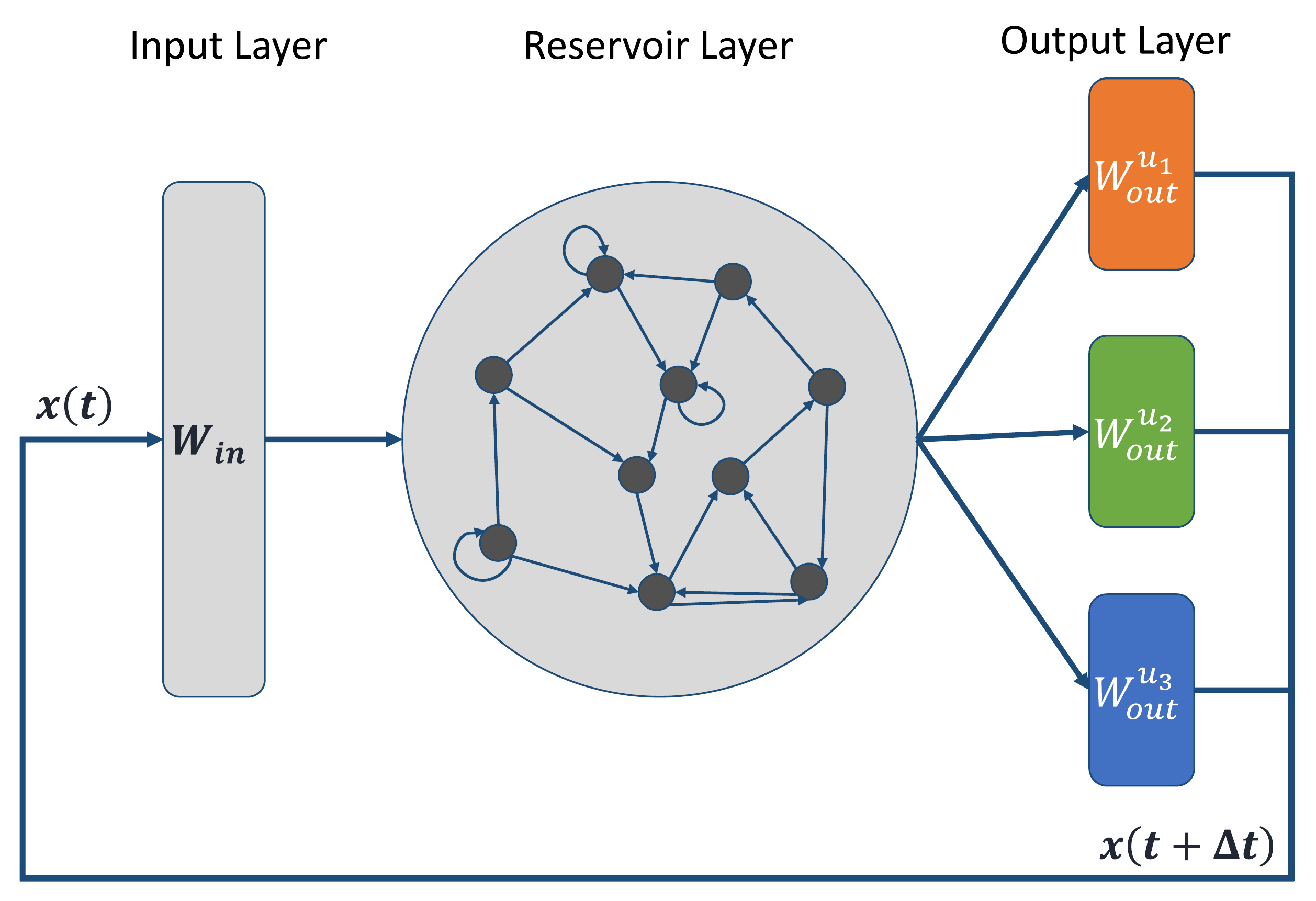}
	\caption{Schematic of an ESN with different readout layers corresponding to different autonomous systems.}
	\label{fig:ESN_ctrl}
\end{figure}

Here, we use an \revKaty{\emph{Echo State Network} (ESN) (or \emph{\revKaty{R}eservoir \revKaty{C}omputer})} \cite{J01} as the surrogate model (cf.\ Figure \ref{fig:ESN_ctrl}), which is a specific type of recurrent neural network (RNN).
The basic idea is to create a large RNN (the reservoir) whose weights are initialized randomly and cannot be trained. Instead, a linear readout layer is responsible for extracting the prediction of the next time step, and it can be trained very efficiently using linear regression.

Since only the output layers $W_\textnormal{out}^{u^j}$ contain information corresponding to the control $u^j$,
we can use a single reservoir for all possible controls. 
Note, that we do not need to explicitly account for the delay coordinates since the ESN is able to capture past dynamics via the feedback of the reservoir state.
Since the control enters linearly -- similar to the Lorenz system -- the interpolated solution is exact. The results in Figure \ref{fig:Results} show that the system can be stabilized at $1.0$.

\section{Discussion of the advantages \rev{and limitations}} \label{sec:DiscussionAdvantages}
\rev{We now briefly discuss the advantages and limitations of our approach to take the effort of reformulating a continuous control problem into a mixed-integer one in order to make easier use of existing surrogate modeling techniques.}

\paragraph*{Advantages:}
For some surrogate modeling techniques, the task of including control inputs can be tackled via state augmentation, i.e., by introducing $\tilde{y} = (y,u)$. However, this is by no means possible for all approaches. Furthermore, state augmentation results in a higher state space dimension and the loss of a low-dimensional subspace (the attractor) the dynamics is restricted to, such that the collection of data becomes more challenging. Our approach presents a remedy to these issues, which we will now discuss in more detail. 

\subsection{Straightforward use of existing autonomous models and the corresponding error bounds}
For some surrogate models -- most importantly projection based approaches such as POD or reduced basis for partial differential equations -- the state space cannot simply be extended, as it is typically a Sobolev space and not $\R^n$. It is thus not clear how to autonomize the dynamics, which becomes even more challenging due to the need to satisfy the boundary conditions \cite{BCB05}. 
Second, if state augmentation is possible, then it may still strongly limit the type of input. For instance, in the case of the Koopman operator, we are restricted to linear inputs \cite{KM18,POR20}.
QuaSiModO avoids these challenges altogether, as we can simply use the autonomous systems in our control framework.

This has another important consequence, as we can readily use error bounds from the autonomous modeling world if they are available (such as \cite{KV02,HO11,SMT+18,LT20,KNH20,NPP+21,ZZ21}).

\subsection{Easier and more efficient data \rev{collection}}
Besides the simplification regarding the data assimiliation task, there may also be an advantage regarding the data requirements, even though this is challenging to assess. Nevertheless, the arguments supporting this assumption are:
\begin{itemize}
    \item the fact that the dynamics of the controlled system is no longer restricted to a low-dimensional subspace, and
    \item the increased state space dimension.
\end{itemize}
To illustrate the latter point, let us consider a setting where we can collect data from arbitrary points in state space.
\begin{example}[Distance-based error estimates]\label{ex:DistBasedErr}
    When thinking about error bounds, some surrogate modeling techniques such as \emph{radial basis functions (RBF)} rely on a \emph{fill distance} $h$, which is the maximum distance between an arbitrary point in state space to the closest point from the training data set\footnote{The related concept of $M$-\emph{boundedness} is used in error estimates for approximations in Reproducing Kernel Hilbert Space (RKHS)  \cite[Theorem 3.3]{CZ07}.}. According to \cite{W2004}, we can assume that $c > 0$ exists such that
    \begin{align}\label{eq:errFillDistance}
        \norm{\Phi_{u} - \Phi^r_{u}} \le \operatorname{exp}\left(\frac{-c}{h}\right) \norm{\Phi_{u}}.
    \end{align}
    Now let us assume that $y(t)\in\R^n$ and $u(t)\in\R^{n_u}$ 
    \rev{and that we want to \revKaty{obtain different} surrogate models \revKaty{which cause nearly the same prediction error, i.e., for which the same error bound holds;}
    a system with an augmented state space $\tilde{y}(t) = (y(t),u(t)) \in \R^{n+n_u}$ \revKaty{(that is hence autonomized)}, and a set of autonomous models for the QuaSiModO approach.
    \revKaty{To make the bounds comparable,}
    assume that in Eq.\ \eqref{eq:errFillDistance}, the constant $c$ and the norms of the right-hand side of the dynamics are identical for the autonomous and the control-augmented systems.
    Then, to obtain training data sets that have the same fill distance $h$,
    }
    we use an equidistant grid of points with $\overline{n}$ points in each spatial dimension and evaluate $\Phi$ once for each point. We thus require $\overline{n}^n$ and $\overline{n}^{n+n_u}$ data points for the $m$ autonomous systems and the autonomized control system, respectively, to achieve the same value for $h$. Assuming that we want to use the control bounds in the QuaSiModO approach, we require $2^{n_u}$ autonomous surrogate models (depending on the application, this number might be further reduced by a suitable choice of $V$). Now, comparing the number of required data points to achieve a comparable bound for the approximation error, we obtain
    \begin{align*}
        \frac{\mbox{\# data points augmented}}{\mbox{\# data points QuaSiModO}}=
        \frac{\overline{n}^{n+n_u}}{2^{n_u} \overline{n}^n} = \left(\frac{\overline{n}}{2}\right)^{n_u},
    \end{align*}
    which grows very quickly for increasing $n_u$. 
\end{example}
\begin{remark}
    Not all surrogate modeling techniques have an exponential growth of the required data with the state space dimension. (In \cite{NPP+21}, for instance, the data requirements have a polynomial scaling with respect to the dimension, however under additional assumptions.) Nevertheless, this is a good indicator that the state space dimension can play a crucial role in the data requirements.
\end{remark}

To further strengthen our point on data efficiency, we perform a numerical experiment as well, where we use \revKaty{ESNs} (cf.\ Section~\ref{subsec:MG}), for the Lorenz system (cf.\ Section~\ref{subsec:Lorenz}) that are trained
\begin{itemize}
    \item as autonomous systems corresponding to the entries in $V$, and
    \item as a single system with an augmented state $\tilde{y} = (y,u)$.
\end{itemize}

\begin{figure}
\centering
\parbox{0.45\textwidth}{
    \centering
	\includegraphics[width=0.45\textwidth]{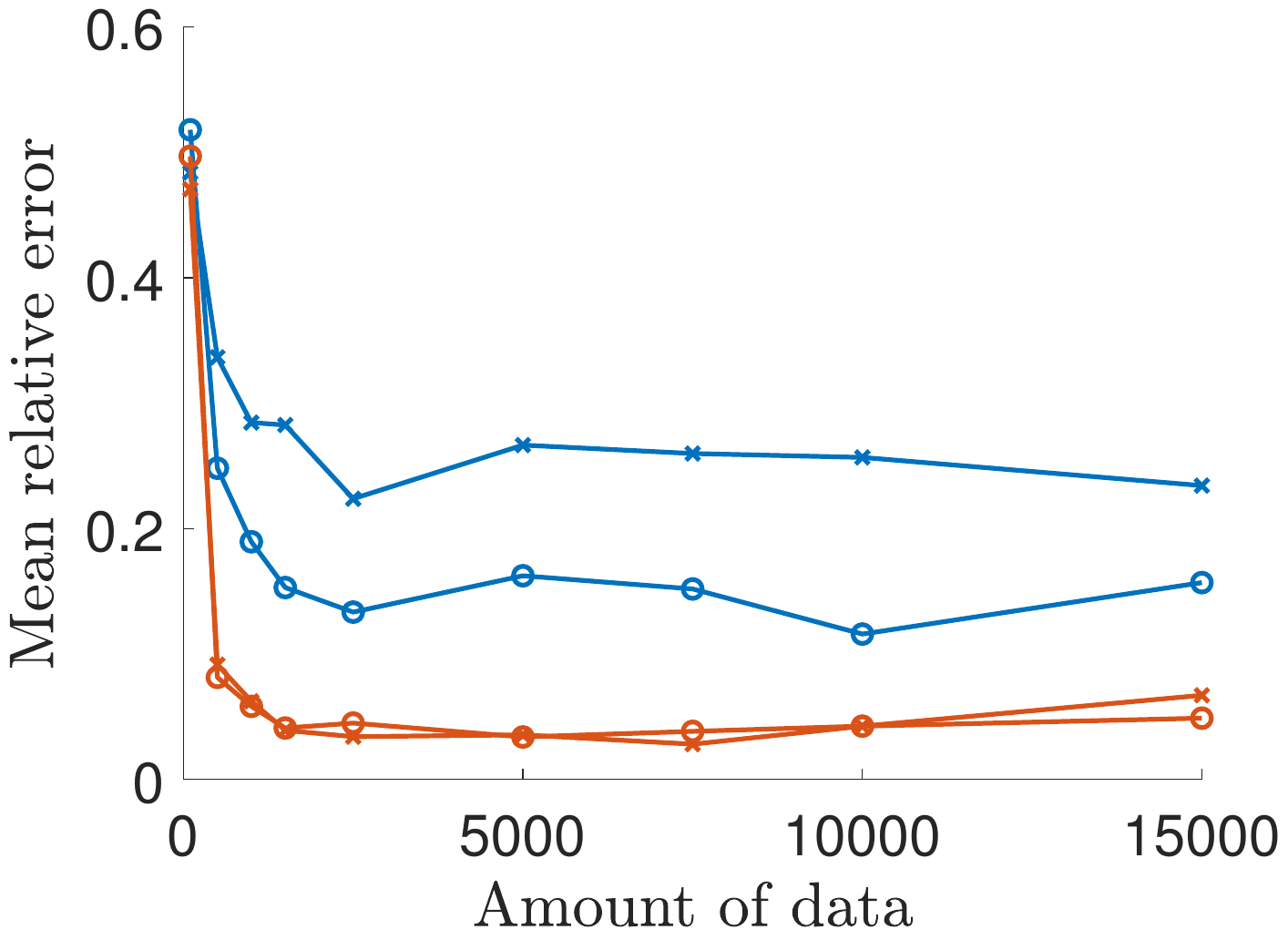}\\
	(a)
}
\hfill
\parbox{0.45\textwidth}{
    \centering
	\includegraphics[width=0.45\textwidth]{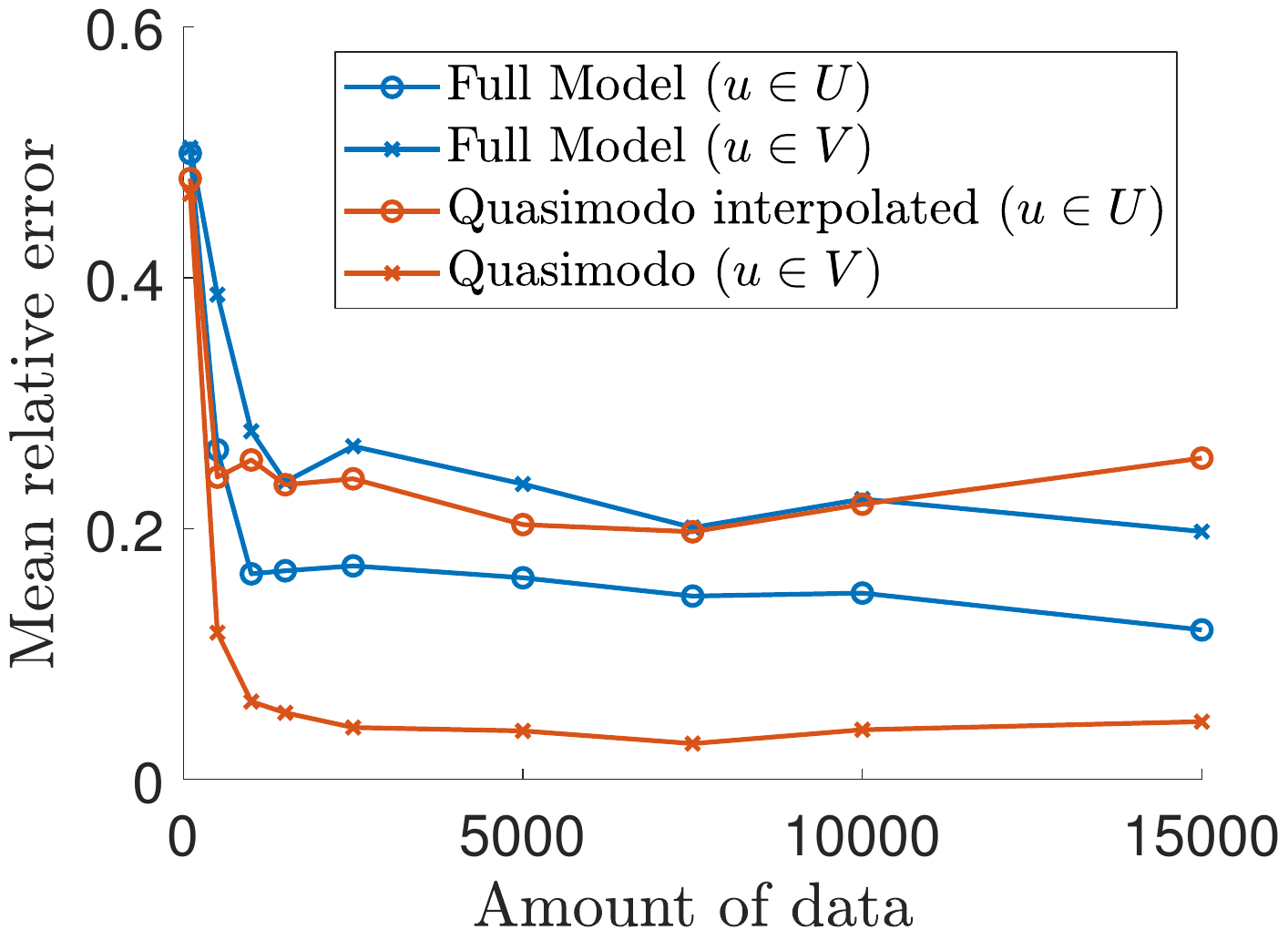}\\
	(b)
}
\caption{Relative $L^2$ prediction error versus the amount of data for (a) the control affine Lorenz system \eqref{eq:Lorenz_linear} and (b) the nonlinear system \eqref{eq:Lorenz_nonlinear}, averaged over $100$ random inputs from $U$ or $V$.} \label{fig:DataReq_ESN}
\end{figure}

\begin{figure}
\centering
\parbox{0.45\textwidth}{
    \centering
	\includegraphics[width=0.45\textwidth]{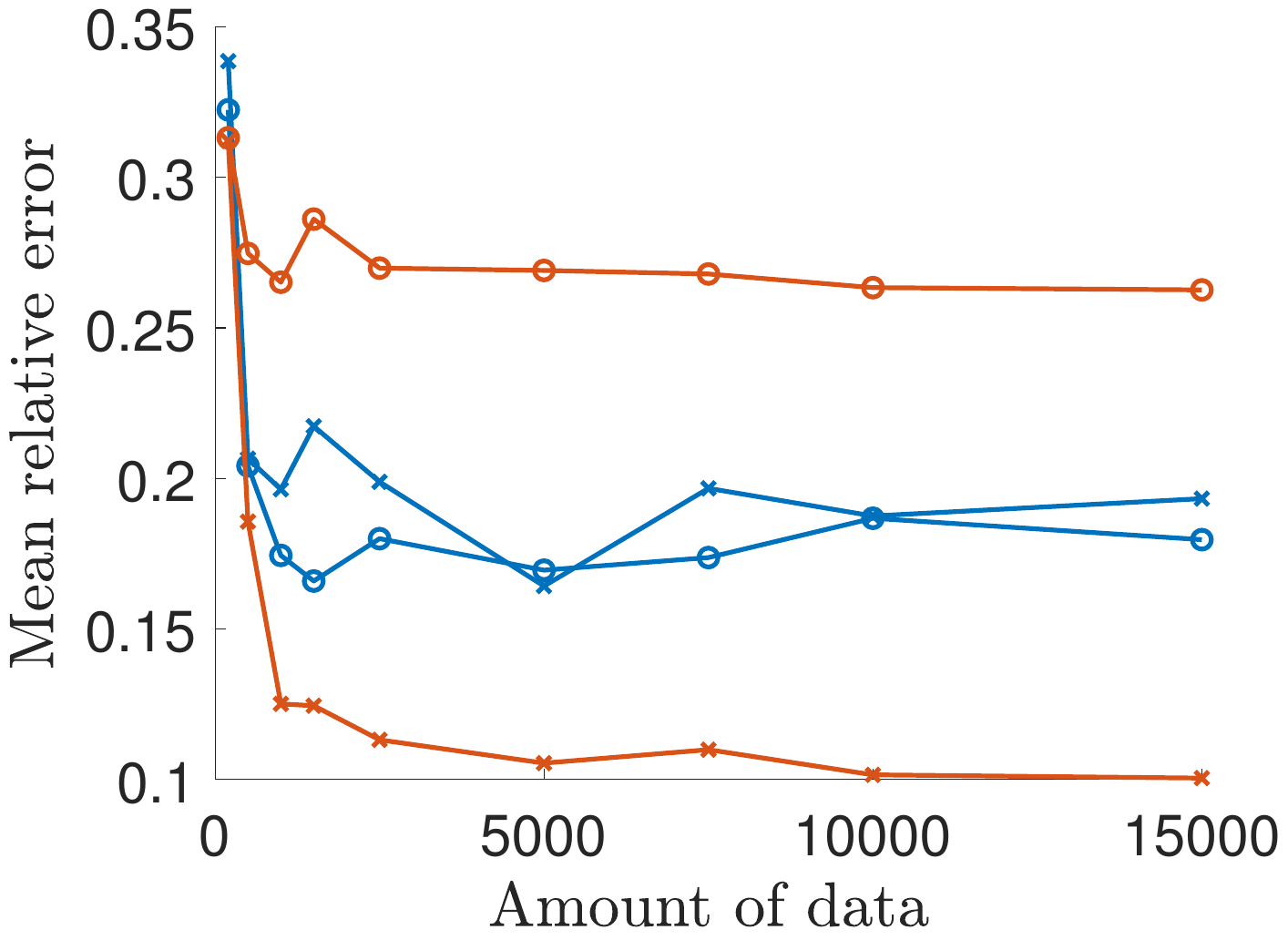}\\
}
\caption{Same as in Figure \ref{fig:DataReq_ESN} (b), only for an LSTM-based surrogate model.}
\label{fig:DataReq_LSTM}
\end{figure}

For a statistically relevant comparison, we repeat the following experiment $100$ times:
\begin{itemize}
    \item Create training data over 750 seconds using random inputs from $V$ or $U$, respectively;
    \item To study the dependency on the size of the training data set, use subsets of different sizes of the collected data (ranging from 100 to \revKaty{15,000} data points);
    \item Use these data sets to train the ESNs for the two approaches; \rev{both approaches are trained with the identical number of samples, but not \revKaty{with the same ones}, as $U\neq V$;}
    \item Compare the relative $L^2$ prediction error over a horizon of 2 seconds, averaged over $100$ simulations with random control sequences from either $U$ or $V$ \revKaty{for both approaches}. 
    For control inputs from $U$, the prediction by the Quasimodo model is realized via linear interpolation of the autonomous systems. 
\end{itemize}
The results are shown in Figure~\ref{fig:DataReq_ESN} for (a) the control affine Lorenz system \eqref{eq:Lorenz_linear} and (b) the modified system with a nonlinear control dependency \eqref{eq:Lorenz_nonlinear}. We observe a clear advantage in the control affine case, not surprisingly regardless of whether the control is selected from $U$ or $V$.
In the nonlinear case, the interpolation error has a negative impact on the QuaSiModO model. However, on the one hand, the performance is still comparable. On the other hand, the sum-up rounding procedure yields controls from $V$ anyway, such that this is the relevant control set for the comparison.
For an LSTM-based surrogate model (cf.\ Figure \ref{fig:DataReq_LSTM}), we observe a very similar behavior.

\begin{figure}
\centering
\parbox{0.45\textwidth}{
    \centering
	\includegraphics[width=0.45\textwidth]{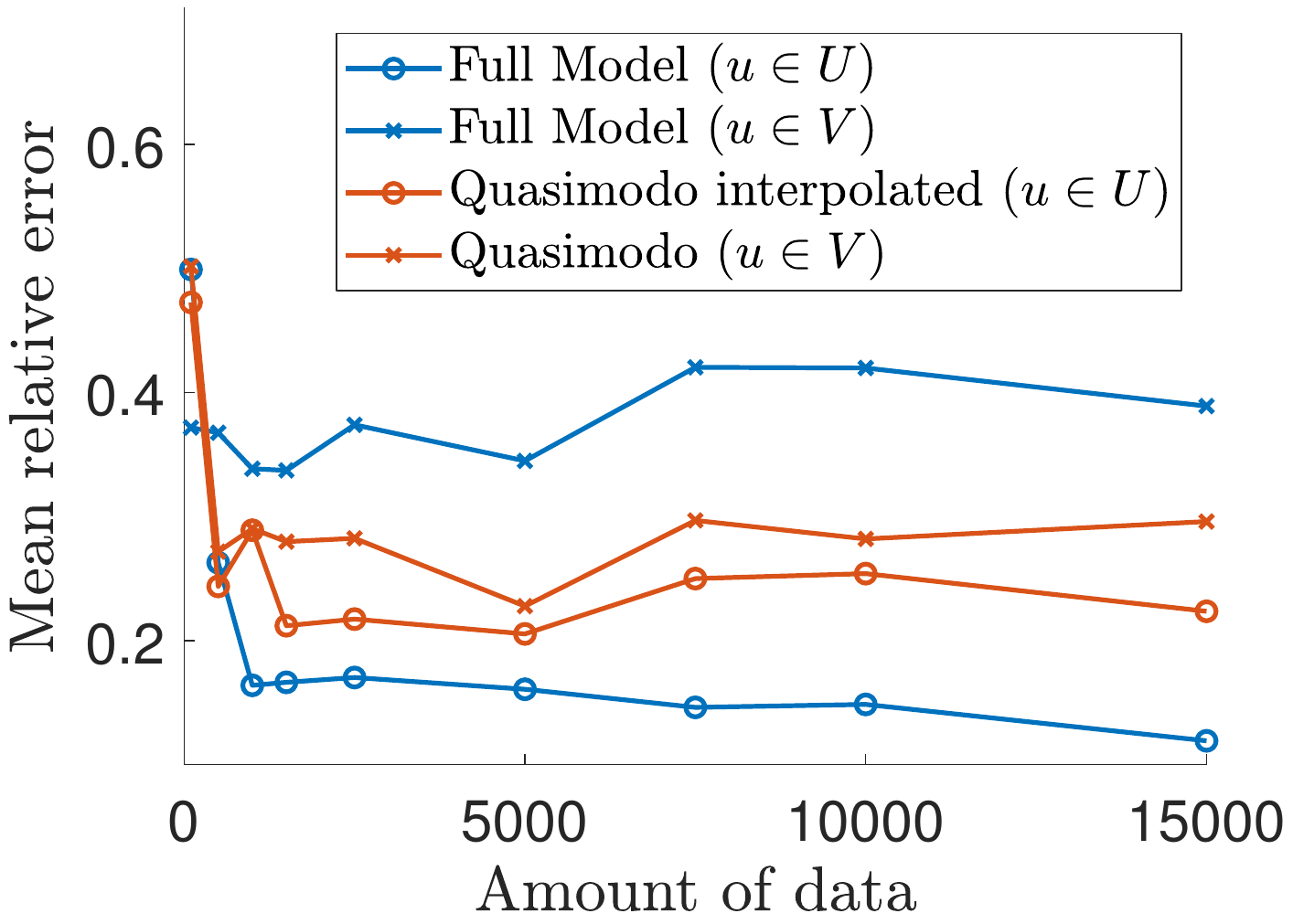}\\
}
\caption{Same as in Figure \ref{fig:DataReq_ESN} (b), but not using the control bounds $U = [0, \pi]$, $V = \{0, 0.8\pi\}$).}
\label{fig:DataReq_ESN_notBounds}
\end{figure}

\rev{
\paragraph*{Limitations:}
Besides the above-mentioned advantages, there are also some limitations -- or at least additional considerations -- that need to be taken into account. Most importantly, we have to mention the task of choosing the finite set $V$. Whereas this is straightforward in the control affine case -- and also in cases with a monotonous dependence on $u$ such as the $\cos(u)$ term in system \eqref{eq:Lorenz_nonlinear} -- this can become more challenging for nonlinear control dependencies. We thus face a challenging decision problem, and if we make a sub-optimal choice, \revKaty{the constant $M_1$ in the error bounds is non-zero} (cf.\ Lemma \ref{lemma:existence_RCCP_solution}) \revKaty{which results in a} reduced control performance.} 
\revKaty{This is caused by the fact that we have a smaller reachable set using only controls from $V$, which becomes also clear if we repeat the analysis from Fig.\ \ref{fig:DataReq_ESN}, but not using the control bounds. The result is demonstrated in Fig.\ \ref{fig:DataReq_ESN_notBounds}, and we see that this has a significant impact on the tracking performance (orange lines).}
\revKaty{However, we would like to mention that many real systems can be modeled in a control affine manner or allow for an informed selection of the control values needed to satisfy the assumptions of Lemma~\ref{lemma:existence_RCCP_solution}.}

\revKaty{If we consider $V$ to be the set of control bounds, we face another, yet related drawback. In this case, we have an}
\rev{
exponential increase in autonomous systems that we need to learn,
i.e., we have $2^{n_u}$ systems that we have to approximate. 
Even though we have demonstrated superiority in data requirements in some cases (cf.\ Example \ref{ex:DistBasedErr}), this may still be computationally very complex. 
\revKaty{For control affine systems, this problem can be prevented if we allow extrapolation.}
This way, we can create $n_u+1$ systems, corresponding to $u=0$ and $n_u$ Euclidean basis vectors for the control space ($u_1 = (1,0,\ldots,0)$, $u_2 = (0,1,0,\ldots,0)$, etc.), which results in a linear complexity with the control dimension. It was already shown that this is a feasible approach for the Koopman operator framework \cite{POR20}.
}

\section{Conclusion}
QuaSiModO is a powerful algorithm for data-driven control of complex systems from many scientific disciplines that does not require knowledge of the underlying dynamics and avoids problem specific modeling efforts. Instead, measurement data corresponding to different fixed inputs can be combined with state-of-the-art surrogate modeling techniques for the prediction of complex dynamical systems in a straightforward manner. 
We demonstrate excellent control performance on a variety of dynamical systems, using different control inputs, observations, and surrogate modeling techniques, thus showing great flexibility and a wide range of possible applications from problems in engineering, biology, or life sciences.
Furthermore, when error bounds are available for the predictive model, these directly translate into error bounds for the corresponding control problem. 
There is a large number of researchers addressing the issue of error bounds for predictive models constructed from data, see, for instance, \cite{SMT+18,LT20}.
Consequently QuaSiModO will directly benefit from these advancements as well as general improvements in data-driven modeling in the future, and will thus continue to play an important role in the construction of data-driven feedback controllers.

\section*{Code} The QuaSiModO toolbox can be obtained freely at \url{https://github.com/SebastianPeitz/QuaSiModO}.

\section*{Acknowledgments}
This research was supported in part by the Priority Programme 1962 of the Deutsche Forschungsgemeinschaft (DFG).

\bibliographystyle{unsrt}
{\footnotesize
\bibliography{bibliography2}
}

\appendix
\section*{Appendix}
\section{Data-driven modeling} \label{app:modeling}

The field of data-driven modeling covers a wide range of algorithms and fields of application. On a very general level, the aim is to use data -- obtained either via numerical simulations or measurements from experiments -- to derive a model that is capable of predicting the future behavior of a system.
The list of surrogate models for dynamical systems is extensive, with additional approaches being presented very regularly. These can be divided into different methodologies:
\begin{itemize}
	\item Projection-based surrogate models \cite{Sir87,HLBR12,BGW15},
	\item Dynamic Mode Decomposition / Koopman operator and generator \cite{Sch10,RMB+09,TRL+14,KGPS18,KM18b,KNP+20,LT20},
	\item (Sparse) regression \cite{CL15,BPK16,RBPK17,SMT+18},
	\item Stochastic modeling approaches \cite{CHK00,BJR11},
	\item Feed-forward neural networks \cite{TF93},
	\item LSTMs \cite{HS97,VBW+18},
	\item Reservoir Computers \cite{JH04,PHG+18}.
\end{itemize}
We have tested several of the above-mentioned approaches within the QuaSiModO framework. These will be described in more detail below.

\subsection{Koopman operator:}
Let $g$ be an autonomous dynamical system defined on the state space $\Y$ (e.g., $\dot{y}(t) = g(y(t))$), and let $ \Phi^t : \Y \to \Y $ with $ \Phi^t(y) = y + \int_0^t g(y(t))\, dt $ be the associated flow map.
Furthermore, let $ f : \Y \rightarrow \R^q $ be a real-valued observable of the system. Then the \emph{Koopman semigroup} of operators $ \{ \Kcal^t \} : \Fcal \to \Fcal $ with $\Fcal = L^{2}(\Y)$, which describes the evolution of the observable $ f $ (see~\cite{Mez13,WKR15,LaMa94,BMM12}), is defined by
\begin{equation*}
(\Kcal^t f)(y) = f(\Phi^t(y)),
\end{equation*}
and the Koopman semigroup is generated by the \emph{Koopman generator} $\Lcal$ \cite{KNP+20,POR20}:
\begin{equation*}
\Lcal f = \lim_{t \to 0} \frac{f\circ \Phi^{t} - f}{t}.
\end{equation*}

The Koopman operator and generator are linear yet infinite-dimensional operators acting on the observable of a system, i.e., on measurements. If we can find a finite-dimensional approximation, then we obtain an entirely data-driven linear system describing the dynamics of the observed quantities $z=f(y)$.
One method to compute such a numerical approximation from data is \emph{Extended Dynamic Mode Decomposition (EDMD)} \cite{WKR15,KKS16} (see \cite{KNP+20} for the generator), which is a generalization of DMD \cite{Sch10,TRL+14}.  
We assume that we have either measurement or simulation data, written in matrix form as
\begin{equation*}
Z = \begin{bmatrix}	z_1 & z_2 & \cdots & z_N \end{bmatrix}
\quad\text{and}\quad
\widetilde{Z} = \begin{bmatrix}	\widetilde{z}_1 & \widetilde{z}_2 & \cdots & \widetilde{z}_N \end{bmatrix},
\end{equation*}
where $z_i = f(y_i)$ and $ \widetilde{z}_i = f(\Phi^{\Delta t}(y_i)) $ with a fixed step size $\Delta t$. 
For a given set of basis functions $ \{ \psi_{1},\,\psi_{2},\,\dots,\,\psi_{k} \} $ (e.g., polynomials, radial basis functions, etc.), the data matrices are embedded into the typically higher-dimensional feature space by
\begin{align*}
\Psi_{Z} = \begin{bmatrix} \psi(z_{1}) & \dots & \psi(z_{m}) \end{bmatrix} \mbox{ and } \Psi_{\widetilde{Z}} =
\begin{bmatrix} \psi(\widetilde{z}_{1}) & \dots & \psi(\widetilde{z}_{m}) \end{bmatrix},
\end{align*}
with $\psi$ being the vector of basis functions.
With these data matrices, we then compute a finite-dimensional approximation of the Koopman operator $ K \in \R^{k \times k} $ by
\begin{equation*}
K^{\top} = \Psi_{\widetilde{Z}} \Psi_Z^+ = \big( \Psi_{\widetilde{Z}} \Psi_Z^{\top} \big) \big(\Psi_Z \Psi_Z^{\top}\big)^{-1},
\end{equation*}
where $^+$ denotes the pseudoinverse. The matrix $K$ now allows us to define a discrete-time update for the observable $z$ which approximates the true dynamics, thus yielding a linear system for the lifted observable $\hat{z} = \psi(z)$:
\[
\hat{z}_{i+1} \approx \Phi^r(\hat{z}_i) = K^\top \hat{z}_i.
\]

\subsection{Proper Orthogonal Decomposition:}
Consider a partial differential equation of the general form 
\begin{equation}
\begin{aligned}
\dot{y}(x,t) &= g(y(x, t)), &(x,t) &\in \Omega \times (t_0, t_e], \\
a(x, t) \frac{\partial y(x, t)}{\partial n} + b(x, t) y(x, t) &= c(x, t), &(x,t) &\in \Gamma \times (t_0, t_e], \\
y(x, t_0) &= y_0(x),&x &\in \Omega,
\end{aligned}
\end{equation}
where $\Omega$ is the spatial domain of interest with boundary $\Gamma = \partial \Omega$ and corresponding outward normal vector $n$. The right-hand side $g$ describes the evolution of the system. For details, the reader is referred to \cite{HPUU09}. Since the state is space dependent, we have to take boundary conditions (BCs) into account, and the coefficients $a(x, t)$, $b(x, t)$ and $c(x, t)$ are given by the problem definition. Note that this covers both Dirichlet as well as Neumann BCs by neglecting one of the terms on the left hand side, respectively.

As a numerical discretization of such PDE systems (e.g., via finite elements) can in general be very expensive to solve, we want to reduce the model dimension by restricting the dynamics to a finite-dimensional linear subspace.
This is realized by representing $y(x,t)$ in terms of basis functions $\{\psi_i(x)\}_{i=1}^{\ell}$:
\begin{align*}
y(x, t) \approx \sum_{i=1}^{\ell} z_i(t) \psi_i(x),
\end{align*}
see \cite{Sir87,HLBR12} 
for details. If we now insert this expansion into the weak formulation of the PDE, we obtain a system of nonlinear ODEs in the coefficients $z$:
\[
\dot{z}(t) = g^r(z(t)) \quad \mbox{and} \quad z_{i+1} = \Phi^r(z_i) = z_i + \int_{t_i}^{t_{i+1}} g^r(z(t))\, \dt.
\]
Note that this requires knowledge of the equations as well as a careful and often tedious treatment of boundary conditions \cite{BCB05} and nonlinearities \cite{CS10}. Moreover, in contrast to the other, non-intrusive modeling techniques, projection-based models predict the full state only.

In order to obtain a small yet informative set of basis functions $\{\psi_i(x)\}_{i=1}^{\ell}$, snapshot data is collected from the (spatially discretized) system and stored in a matrix $\hat{Y}=[\hat{y}(t_0), \ldots, \hat{y}(t_N)]$, where the hat notation denotes the finite-dimensional spatial discretization of $y$. The singular value decomposition of $\hat{Y}$ then yields the orthonormal basis (given by the singular vectors) with the smallest $L^2$ projection error (equal to the sum of the neglected singular values). This procedure is known as \emph{Proper Orthogonal Decomposition} (POD) and has been successfully applied to a wide range of nonlinear systems over the years, fluid dynamics being the most prominent example \cite{HLBR12}. 

Quite a number of extensions have been presented for the use of POD in PDE-constrained control, see, for instance, \cite{KV99} 
for the Burgers equation, which we use as one of our examples.

\subsection{Reservoir Computing / Echo State Networks:}\label{app:RC}

Reservoir Computers (RC) -- also referred to as Echo State Networks (ESN) -- have become very popular for time series prediction. This is mainly due to the straightforward and fast training process, which only consists of solving a linear system. Further details on RC and ESN can be found in the survey \cite{LJ09} or in \cite{JH04,J02,J01}.

The basic idea is to create a large recurrent neural network whose weights are initialized randomly and cannot be trained. This is called the reservoir. 
A linear output layer is then added, which can be trained efficiently via linear regression. The standard equations of an ESN are given by
\begin{align*}
r(k + 1) &= \sigma(W^{\mathsf{in}} i(k)  + W^{\mathsf{res}} r(k) + W^{\mathsf{fb}} o(k)),\\
o(k + 1) &= W^{\mathsf{out}}r(k + 1), 
\end{align*}
where $\sigma$ is a nonlinear activation function, e.g., $\sigma(x) = \tanh(x)$, $W^{\mathsf{in}}$, $W^{\mathsf{res}}$ and $W^{\mathsf{fb}}$ are randomly generated (sparse) matrices and $W^{\mathsf{out}}$ is the trainable output matrix. 
The reservoir state $r(k+1)$ (for time step $t_{k+1}$) is computed based on a time-dependent input $i(k)$, the previous reservoir state $r(k)$ and the previous output $o(k)$. In the context of time series prediction the output is the (approximated) state of the system, i.e., $o(k) = y_k$. Since there is no additional time-dependent input in this setting, the term $W^{\mathsf{in}} i(k)$ is omitted. 

To train the model, we take -- similar to the Koopman operator approach -- observed data (from measurements or simulations) to train our model, i.e., $\Phi^t : \Y \to \Y$ is the flow map which describes the dynamics of the system and $ f : \Y \rightarrow \R^q $ is a real-valued observable of the system. As the reservoir has a ``memory'' due to the feedback of the reservoir state, the data to train an ESN has to stem from a single time series:
\begin{align}
Z = \begin{bmatrix}	z_{t_{0}} & z_{t_{1}} & \cdots & z_{t_N}   \end{bmatrix} \text{ with } z_{t_k} = f(\Phi^{t_k}(y_0)).
\end{align}
The first time steps are usually used to initialize the reservoir and are not used to train the linear output layer. Furthermore, in the feedback loop the output $o(k)$ is replaced by the true system state from the training data $z_{t_k}$.

There are many publications on the prediction of chaotic dynamical systems with ESNs, see, e.g., \cite{LHW12,MW17,CHS19,PLH17,LHO18,J01}. Furthermore, some papers discuss the use of ESN in an MPC context \cite{ATF19, JCA18}. Therein, the control input $u_k$ serves as the input $i(k)$. Due to this, the input space dimension increases, and therefore, a larger number of neurons is required in the reservoir to maintain a sufficiently accurate prediction. Moreover, the training will become harder. 

%

\subsection{Long-short-term memory (LSTM) neural networks:}

LSTM is a specific architecture for neural networks, more precisely for recurrent neural networks, which is specifically tailored to sequential data \cite{HS97}, 
e.g., time series prediction. In \cite{VBW+18}, 
the authors successfully applied the LSTM-approach to forecasting chaotic systems in a reduced order space. 
Here, we use the standard \emph{tensorflow} implementation \cite{TF} which coincides with the formulation in \cite{HS97}. 

Consider the flow map $\Phi^t : \Y \to \Y$ which describes the dynamics of the system.
Furthermore, let $ f : \Y \rightarrow \R^q $ be a real-valued observable. We assume that we have training data where one data point is of the following form:
\begin{align}
Z_{\mathsf{in}} = \begin{bmatrix}	z_{t_{k-d}} & z_{t_{k-d+1}} & \cdots & z_{t_k} \end{bmatrix}
\quad\text{and}\quad
Z_{\mathsf{out}} = z_{t_{k+1}}
\end{align}
with $z_{t_k} = f(\Phi^{t_k}(y_0))$ for a given $y_0 \in \Y$, i.e., the LSTM gets the time series $Z_{\mathsf{in}}$ as a sequential input and should predict the behavior of the dynamical system for one time step into the future ($Z_{\mathsf{out}}$ ).
Note that in the presented control framework, for training an LSTM model corresponding to control $u^j$, the delayed time series can be produced by different control inputs, but only the control $u_k$ which maps $z_{t_k}$ to $z_{t_{k+1}}$ has to be equal to $u^j$.

\section{Detailed description of the numerical examples}\label{app:numerics}
In this section, we give a detailed description of the numerical setup of the examples presented in Figure \ref{fig:Results}. In all examples, we solve Problem \eqref{eq:OCP-IV} with a tracking objective:
\begin{equation} \tag{$\widehat{\text{IV}}$} \label{eq:OCP-IVr}
\begin{aligned}
\min_{\alpha\in ([0,1]^m)^p} &\sum_{i=0}^{p-1} (z_{i+1}-z_{i+1}^{\mathsf{ref}})^\top Q (z_{i+1}-z_{i+1}^{\mathsf{ref}})\\ 
\mbox{s.t.} \quad z_{i+1} = \Phi^r(z_i, \alpha_i) &= \sum_{j=1}^m \alpha_{i,j} \Phi^r_{u^j}(z_i) \quad \mbox{and} \quad
\sum_{j=1}^m \alpha_{i,j} = 1, \qquad i = 0,1,2,\ldots,
\end{aligned}
\end{equation}
where the quadratic, positive semidefinite matrix $Q$ is problem-specific.

\subsection{Control of the Lorenz system using the Koopman operator:}\label{subsec:Lorenz}
To complete the description in the results section, a detailed setting of the numerical setup is described in Table~\ref{tab:Param_Lorenz_EDMD}.
\begin{table}[b!]
	\footnotesize
	\centering
	\caption{Lorenz system with EDMD surrogate model.}
	\label{tab:Param_Lorenz_EDMD}
	\begin{tabular}{lll}
		\toprule
		& Parameter & Value\\
		\midrule
		System parameters & $(\sigma, \rho, \beta)$ & $(10, 28, \frac{8}{3})$ \\
		\midrule
		Quantization & $U$ & $[-50, 50]$ or $[0, \pi]$\\
		& $V$ &  $\{-50, 50\}$ or $\{0, \pi\}$\\
		& $m$ & $2$ \\
		\midrule
		Training data & $\Delta t$ & $0.0005$ \\
		& $T_{\mathsf{train}}$ & $100.0$ \\
		& \# trajectories & $1$ \\
		& Input & Piecewise constant, random $u_i\in V$\\
		\midrule
		Surrogate model & $\Delta t$ & $0.05$ \\
		& $\psi$ & Monomials up to degree $3$\\ 
		& observable & $z=f(y)=y$ \\ 
		\midrule
		MPC & $T_{\mathsf{MPC}}$ & $20.0$ \\
		& $p$ & $3$\\ 
		& $Q$ & $\operatorname{diag}(0,1,0)$\\
		& $z_i^{\mathsf{ref}}$ & $z_i^{\mathsf{ref}} = 1.5 \cdot \sin(4\pi \cdot t_i / T_{\mathsf{MPC}})$\\ 
		& $\Delta t_{\mathsf{SUR}}$ & $0.0005$\\
		\bottomrule
	\end{tabular}
\end{table}

\subsection{Control of the Burgers equation using Proper Orthogonal Decomposition:}

Our second example is the one-dimensional viscous Burgers equation:
\begin{align*}
\dot{y}(x,t) - \frac{1}{Re} \Delta y(x,t) + y(x,t) \nabla y(x,t) &= \sum_{j=1}^5 v^j(t) \chi^j(x), \quad &(x,t) \in [0,L] \times (t_0, t_e],\\
y(0, t) = y(L, t) &= 0, \quad &t \in (t_0, t_e], \\
y(x,t_0) &= \begin{cases} 1, & x\in (0, \frac{L}{2}], \\ 0, & x \in (\frac{L}{2},L),  \end{cases}
\end{align*}
which was also studied by Kunisch and Volkwein in their seminal work on POD-based control of PDEs \cite{KV99}. 
We consider a domain of length $L=1$ and a distributed control that is realized via indicator functions $\chi_i$ with disjoint support:
\[
\chi^j(x) = \begin{cases} 1, & \frac{(j-1)L}{5} < x \leq \frac{jL}{5} \\ 0, & \mbox{else}  \end{cases} \quad \Rightarrow \quad \sum_{j=1}^5 \chi^j(x) = 1 ~\mbox{for}~x\in(0,L],
\]
with the aim to stabilize the system at $y^{\mathsf{ref}}(\cdot,t) = 0$. 
For the quantization of the control, we use a difference-star-like set with $u^1 = (0,0,0,0,0)$ and then two additional points per component, where the minimal and maximal value are taken, respectively. This results in a total of 11 autonomous systems.
As we use POD as the surrogate model, we have to choose $z=y$. The detailed setup is described in Table \ref{tab:Param_Burgers_POD}.

\begin{table}[h!]
	\centering
	\footnotesize
	\caption{1D Burgers equation with POD surrogate model.}
	\label{tab:Param_Burgers_POD}
	\begin{tabular}{lll}
		\toprule
		& Parameter & Value\\
		\midrule
		System parameters & $Re$ & 100 \\
		\midrule
		Quantization & $U$ & $[-1, 1]^5$ \\
		& $V$ &  $\left\{\begin{pmatrix} -1 \\ 0\\ \vdots \\ 0 \end{pmatrix}, \begin{pmatrix} 1 \\ 0\\ \vdots \\ 0 \end{pmatrix}, \ldots, \begin{pmatrix} 0 \\ \vdots \\ 0 \\ 1 \end{pmatrix} \right\}$ \\
		& $m$ & 11 \\
		\midrule
		Training data & $\Delta t$ & 0.005 \\
		& $T_{\mathsf{train}}$ & 50 \\
		& \# trajectories & 1 \\
		& Input & Piecewise constant, random $u_i\in V$\\
		\midrule
		Surrogate model & $\Delta t$ & 0.025 \\
		& Basis size $\ell$ & 12 \\
		& observable & $z=f(y)=y$ \\ 
		\midrule
		MPC & $T_{\mathsf{MPC}}$ & 5 \\
		& $p$ & 5 \\ 
		& $Q$ & $\mathsf{Id}$\\
		& $z_i^{\mathsf{ref}}$ & $0$\\ 
		\bottomrule
	\end{tabular}
\end{table}

\subsection{Control of the flow around a cylinder using an LSTM neural network:}
Next, we want to control the flow around a cylinder -- governed by the two-dimensional incompressible Navier--Stokes equations -- which exhibits periodic vortex shedding at $Re=100$ \cite{NAM+03}.
Our aim is not to control the flow field, but the forces acting on the cylinder (the lift $C_L$ and the drag $C_D$) without any knowledge of the flow field itself.

As the surrogate model we use an LSTM neural network.
The control aim is to force the lift to follow a predefined trajectory. 
This was already done in \cite{BPB+20} 
but there the control input was continuous and a different RNN architecture was used to build the surrogate model. For detailed parameters of the setting see Table~\ref{tab:Param_Cylinder_LSTM}.

\begin{table}[h!]
	\footnotesize
	\centering
	\caption{2D Cylinder flow with LSTM surrogate model.}
	\label{tab:Param_Cylinder_LSTM}
	\begin{tabular}{lll}
		\toprule
		& Parameter & Value\\
		\midrule
		System parameters & $Re$ & 100 \\
		\midrule
		Quantization & $U$ & $[-5, 5]$ \\
		& $V$ &  $\{-5, 0, 5\}$ \\
		& $m$ & 3 \\
		\midrule
		Training data & $\Delta t$ & 0.05 \\
		& $T_{\mathsf{train}}$ & (2 *) 500 \\
		& \# trajectories & 2 \\
		& Input & Piecewise constant, random $u_i\in V$\\
		\midrule
		Surrogate model & $\Delta t$ & 0.1 \\
		& $\text{neurons per LSTM-cell}$ & 500 \\
		& observable & $z = f(y) = (C_D,C_L)$\\ 
		& delay coordinates & 15\\ 
		\midrule
		MPC & $T_{\mathsf{MPC}}$ & 20.0 \\
		& $p$ & 5 \\ 
		& $Q$ & $\begin{bmatrix}
		0 & 0\\
		0 & 1 
		\end{bmatrix}$ \vspace*{0.1cm}\\
		& $z_i^{\mathsf{ref}}$ & $\sin(\frac{t_i}{2})$\\ 
		& $\Delta t_{\mathsf{SUR}}$ & 0.05\\
		\bottomrule
	\end{tabular}
\end{table}

\subsection{Control of the Fluidic Pinball using the Koopman operator:}
As before, our aim is not to control the flow field, but the forces acting on the (this time three) cylinders without any knowledge of the flow field itself.

As the surrogate model we again use the Koopman operator, and the control aim is to force the lift to follow a predefined trajectory abs before. For detailed parameters of the setting see Table~\ref{tab:Param_Cylinder_LSTM}.

\begin{table}[h!]
	\footnotesize
	\centering
	\caption{2D Cylinder flow with LSTM surrogate model.}
	\label{tab:Param_Cylinder_LSTM}
	\begin{tabular}{lll}
		\toprule
		& Parameter & Value\\
		\midrule
		System parameters & $Re$ & 200 \\
		\midrule
		Quantization & $U$ & $[-2, 2]^3$ \\
		& $V$ &  $\{-2, 2\}^3$ \\
		& $m$ & 8 \\
		\midrule
		Training data & $\Delta t$ & 0.05 \\
		& $T_{\mathsf{train}}$ & 100 \\
		& \# trajectories & 1 \\
		& Input & Piecewise constant, random $u_i\in V$\\
		\midrule
		Surrogate model & $\Delta t$ & $0.25$ \\
		& $\psi$ & Monomials up to degree $1$\\ 
		& observable & $z = f(y) = (C_{D,1},C_{L,1},C_{D,2},C_{L,2},C_{D,3},C_{L,3})$\\
		& delay coordinates & 1\\ 
		\midrule
		MPC & $T_{\mathsf{MPC}}$ & 60.0 \\
		& $p$ & 5 \\ 
		& $Q$ & $\operatorname{diag}(0,1,0,1,0,1)$\\
		& $z_i^{\mathsf{ref}}$ & piecewise constant\\ 
		\bottomrule
	\end{tabular}
\end{table}

\subsection{Control of the Mackey-Glass delay differential equation for blood cell reproduction using Reservoir Computing (ESN):}
The last example is the control of the Mackey-Glass equation which is a delay differential equation modeling blood cell reproduction \cite{MG77,GM79}: 
\begin{equation}
\begin{aligned}
\dot{y}(t) = \beta \frac{y(t -\tau)}{1 + y(t -\tau)^\eta} - \gamma y(t) + u(t) \quad \text{ with } \beta, \gamma, \eta > 0.
\end{aligned}
\end{equation}
The uncontrolled system (i.e., $u(t) = 0$) was studied for different parameters and chaotic behavior was proven for certain parameter values, for instance $\beta = 2$, $\gamma = 1$ and $\eta = 9.65$ \cite{GM79}.
Since this system has become a benchmark for predicting chaotic delay systems, there are several studies on data-based models, and many different methods were used such as local linear approximation \cite{FS1987}, genetic algorithms \cite{FTK2006} and different forms of neural networks \cite{MBS1993, L_pez_Caraballo_2016,FLPK2019,ZLXZ2014}. Another approach is via Echo State Networks (ESN), which were studied in \cite{J01, GSS2015, HL2006}.

\begin{table}[t]
	\footnotesize
	\centering
	\caption{Mackey-Glass DDE with ESN surrogate model.}
	\label{tab:Param_MackeyGlass_ESN}
	\begin{tabular}{lll}
		\toprule
		& Parameter & Value\\
		\midrule
		System paramters & $(\beta, \gamma, \eta, \tau)$ & $(2, 1, 9.65, 2)$ \\
		\midrule
		Quantization & $U$ & $[-0.2, 1.0]$ \\
		& $V$ &  $\{-0.2,0.0, 1.0\}$ \\
		& $m$ &  $3$ \\
		\midrule
		Training data & $\Delta t$ & $0.05$ \\
		& $T_{\mathsf{train}}$ & $1000.0$ \\
		& \# trajectories & 1 \\
		& Input & Piecewise constant, random $u_i\in V$\\
		\midrule
		Surrogate model & $\Delta t$ & 0.25 \\ & $\text{size residuum}$ & 200 \\
		& $\text{spectral radius}$ ($W^{\mathsf{res}}$) & $0.75$\\ 
		& $\text{sparsity}$ ($W^{\mathsf{res}}$)  & $0.9$\\ 
		& $\sigma$ & $0.99$\\ 
		& $\beta$ & $0.0001$\\ 
		& observable & $z= f(y) = y$ \\ 
		\midrule
		MPC & $T_{\mathsf{MPC}}$ & $20.0$ \\
		& $p$ & 5\\ 
		& $Q$ & $1$\\
		& $z_i^{\mathsf{ref}}$ & $1.0$\\ 
		\bottomrule
	\end{tabular}
\end{table}

The additive term to control the system represents -- in the context of blood reproduction -- an increase in the number of blood cells caused by, e.g., a transfusion. The same control problem was already studied in \cite{SDL2011,Kiss2017}, where the authors derived different feedback laws in order to stabilize the systems.

Here, we use a ESN as surrogate model. As mentioned in Appendix~\ref{app:RC} we can use a single reservoir for all possible controls, i.e., only the output layer $W^{\mathsf{out}}$ differs. For the detailed setting see Table~\ref{tab:Param_MackeyGlass_ESN}. Note, that we do not need delay coordinates although the state depends on $y(t)$ and $y(t-\tau)$ since the ESN is able to capture the past dynamics by the feedback of the reservoir state.     
Since the control enters linearly -- similar to the Lorenz system -- the interpolated solution is exact ($\Ev= 0$). The results in Figure \ref{fig:Results} show that the system can be stabilized at $1.0$. 

\end{document}